\documentclass[11pt,reqno]{amsart}
\usepackage{amsmath,amsfonts,amsthm,amsopn,amssymb,extarrows}
\usepackage{cite,marginnote}
\usepackage{color,enumitem,graphicx}
\usepackage[colorlinks=true,urlcolor=blue,
citecolor=red,linkcolor=blue,linktocpage,pdfpagelabels,
bookmarksnumbered,bookmarksopen]{hyperref}
\usepackage[english]{babel}

\usepackage[left=2.9cm,right=2.9cm,top=2.8cm,bottom=2.8cm]{geometry}
\usepackage[hyperpageref]{backref}

\numberwithin{equation}{section}

\makeindex
\makeindex
\newtheorem{theorem}{Theorem}[section]
\newtheorem{definition}[theorem]{Definition}
\newtheorem{lemma}[theorem]{Lemma}
\newtheorem{corollary}[theorem]{Corollary}
\newtheorem{proposition}[theorem]{Proposition}
\newtheorem{remark}[theorem]{Remark}

\newcommand{\ud}{\mathrm{d}}

\newcommand{\s}{\section}

\newcommand{\R}{\mathbb R}

\newcommand{\lab}{\label}
\newcommand{\bt}{\begin{theorem}}
	\newcommand{\et}{\end{theorem}}
\newcommand{\bl}{\begin{lemma}}
	\newcommand{\el}{\end{lemma}}
\newcommand{\bd}{\begin{definition}}
	\newcommand{\ed}{\end{definition}}
\newcommand{\bc}{\begin{corollary}}
	\newcommand{\ec}{\end{corollary}}
\newcommand{\bp}{\begin{proof}}
	\newcommand{\ep}{\end{proof}}
\newcommand{\bx}{\begin{example}}
	\newcommand{\ex}{\end{example}}
\newcommand{\bi}{\begin{exercise}}
	\newcommand{\ei}{\end{exercise}}
\newcommand{\bo}{\begin{proposition}}
	\newcommand{\eo}{\end{proposition}}
\newcommand{\br}{\begin{remark}}
	\newcommand{\er}{\end{remark}}
\newcommand{\beq}{\begin{equation}}
	\newcommand{\eeq}{\end{equation}}
\newcommand{\ba}{\begin{align}}
	\newcommand{\ea}{\end{align}}
\newcommand{\bn}{\begin{enumerate}}
	\newcommand{\en}{\end{enumerate}}
\newcommand{\bg}{\begin{align*}}
	\newcommand{\bcs}{\begin{cases}}
		\newcommand{\ecs}{\end{cases}}

	\newcommand{\bean}{\begin{eqnarray*}}
		\newcommand{\eean}{\end{eqnarray*}}

	\def\C{\mathbb{C}}
	\def\N{\mathbb{N}}

	\def\R{\mathbb{R}}

	\def\bd{\mathrm{bd}\,}

	\newcommand{\cf}{{\mathcal F}}

\title[Normalized solutions of coupled Sobolev critical]{\small
	 Normalized solutions of coupled Sobolev critical Schr\"odinger equations with mass subcritical couplings}

	\author[J.~J.~Zhang]{Jianjun Zhang}
	\author[X.~X.~Zhong]{Xuexiu Zhong}
	\author[J.~F.~Zhou]{Jinfang Zhou}

\address[J.~J.~Zhang]{\newline\indent College of Mathematics and Statistics
\newline\indent
Chongqing Jiaotong University
\newline\indent
Xuefu, Nan'an, 400074, Chongqing, PR China}
\email{\href{mailto:zhangjianjun09@tsinghua.org.cn}{zhangjianjun09@tsinghua.org.cn}}

	\address[X.~X.~Zhong]{\newline\indent South China Research Center for Applied Mathematics and Interdisciplinary Studies
		\newline\indent
		South China Normal University
		\newline\indent
		Guangzhou 510631, P. R. China}
	 \email{\href{mailto:zhongxuexiu1989@163.com}{zhongxuexiu1989@163.com}}
	
	\address[J.~F.~Zhou]{\newline\indent School of Mathematical Sciences
		\newline\indent
		South China Normal University
		\newline\indent
		Guangzhou 510631, P. R. China}
	 \email{\href{mailto:z18870300649@163.com}{z18870300649@163.com}}
	
\begin{document}
		
\begin{abstract}
We are concerned with qualitative properties of positive solutions to the following coupled Sobolev critical Schr\"odinger equations
$$
\begin{cases}
-\Delta u+\lambda_1 u=\mu_1|u|^{2^*-2}u+\nu\alpha |u|^{\alpha-2}|v|^{\beta}u ~\hbox{in}~ \R^N,\\
-\Delta v+\lambda_2 v=\mu_2|v|^{2^*-2}v+\nu\beta |u|^{\alpha}|v|^{\beta-2}v ~\hbox{in}~ \R^N
\end{cases}
$$
subject to the mass constraints $\int_{\mathbb{R}^N}|u|^2 \ud x=a^2$ and $\int_{\mathbb{R}^N}|v|^2 \ud x=b^2$, where, $a>0,\,b>0,\,N=3,4$ and $2^*:=\frac{2N}{N-2}$ is the Sobolev critical exponent.  The main purpose of this paper is focused on the mass mixed case, i. e., $ \alpha>1,\beta>1,\alpha+\beta<2+\frac{4}{N}$.  For some suitable small $\nu>0$, we show that the above system  admits two positive solutions, one of which is a local minimizer, and another one is a mountain pass solution. Moreover, as $\nu\to0^+$, asymptotic behaviors of solutions are also considered. Our result gives an affirmative answer to a Soave's type open problem raised by Bartsch {\it et al.} (Calc. Var. Partial Differential Equations 62(1), Paper No. 9, 34, 2023).

\begin{flushleft}
\end{flushleft}
\begin{flushleft}
\end{flushleft}

\end{abstract}
\maketitle

\s{Introduction and main results}\label{Introduction}
The time-dependent system of coupled nonlinear Schr\"odinger equations
\beq\label{eq:230905-1}
\begin{cases}
	-i\partial_t\Psi_1=\Delta \Psi_1+\mu_1|\Psi_1|^{p-2}\Psi_1+\nu\alpha|\Psi_1|^{\alpha-2}|\Psi_2|^{\beta}\Psi_1,\\
	-i\partial_t\Psi_2=\Delta \Psi_2+\mu_2|\Psi_2|^{q-2}\Psi_2+\nu\beta|\Psi_1|^{\alpha}|\Psi_2|^{\beta-2}\Psi_2,\\
	\Psi_j=\Psi_j (x,t)\in\C, (x,t)\in \R^N\times\R, j=1,2,
\end{cases}
\eeq
 has received widespread attention from scholars due to its rich physical background.
It is derived from the mean field model of binary mixtures of Bose-Einstein condensates or binary gases of  fermion atoms degenerate quantum state (Bose-Fermi mixtures, Fermi-Fermi mixtures), see e.g. \cite{Adhikari2007,Bagnato2015,Esry1997,Malomed2008}.
An important, and of course well known feature of system \eqref{eq:230905-1}, is the conservation of masses: the $L^2$-norms $|\Psi_1(\cdot,t)|_2, |\Psi_2(\cdot,t)|_2$ of solutions are independent of $t\in\R$. These norms have a clear physical meaning. In the aforementioned contexts, they represent the number of particles of each component in Bose-Einstein condensates, or the power supply in the nonlinear optics framework.

A solitary wave solution of \eqref{eq:230905-1} with the form $\Psi_1(x,t)=e^{i\lambda_1 t}u_1(x)$ and $\Psi_2(x,t)=e^{i\lambda_2 t}u_2(x)$ leads to the following elliptic system
\beq\label{eq:230905-2}
\begin{cases}
	-\Delta u+\lambda_1 u=\mu_1|u|^{p-2}u+\nu\alpha |u|^{\alpha-2}|v|^{\beta}u, ~\hbox{in}~ \R^N,\\
	-\Delta v+\lambda_2 v=\mu_2|v|^{q-2}v+\nu\beta |u|^{\alpha}|v|^{\beta-2}v, ~\hbox{in}~ \R^N.
\end{cases}
\eeq
A solution satisfying the prescribed $L^2$-norm that $\int_{\mathbb{R}^N}|u|^2 \ud x=\|\Psi_1(\cdot,0)\|_2^2$ and $\int_{\mathbb{R}^N}|v|^2 \ud x=\|\Psi_2(\cdot,0)\|_2^2$
is called a normalized solution in the literature. Due to the physical relevance of normalized solutions, the study of normalized solutions for a nonlinear Schr\"odinger equation or system is a hot topic in recent years. A natural approach to finding solutions of \eqref{eq:230905-2} satisfying the normalization constraints
\beq\label{eq:230814-2}
\int_{\R^N}u^2\ud x=a^2,\quad \int_{\R^N}v^2\ud x=b^2
\eeq
consists in finding critical points $(u,v)\in H^1(\R^N, \R^2)$ of the energy
$$I(u,v)=\frac{1}{2}\int_{\R^N}(|\nabla u|^2+|\nabla v|^2)\ud x-\frac{\mu_1}{p}\int_{\R^N}|u|^p\ud x-\frac{\mu_2}{q}|v|^q\ud x-\nu\int_{\R^N}|u|^{\alpha}|v|^{\beta}\ud x$$
under the constraints \eqref{eq:230814-2}. In this kind of approach, $\lambda_1,\lambda_2$ are unknown parameters which will come out as Lagrange multipliers.

Comparing to the \textit{fixed frequency} problem, the study of normalized solution problems
 possesses the following technical difficulties when dealing with it in the variational framework:
 \begin{itemize}
\item[(1)] One can not use the usual Nehari manifold method since the frequency is unknown.
\item[(2)] The existence of bounded Palais-Smale sequences requires new arguments.
\item[(3)]  The Lagrange multipliers have to be controlled.
\item[(3)] The embedding $H^1(\R^N)\hookrightarrow L^2(\R^N)$ is not compact. Even if we prefer to work in $H_{rad}^{1}(\R^N)$, the embedding $H^{1}_{rad}(\R^N)\hookrightarrow L^2(\R^N)$ is also not compact. For the {\it fixed frequency problem}, usually a nontrivial weak limit is also a solution. However, for the {\it fixed mass problem}, even the weak limit is nontrivial, the constraint condition may be not satisfied.
\item[(4)] The number $\bar{p}=2+\frac{4}{N}$, called the mass critical exponent, affects the geometry of the functional heavily.
\end{itemize}

For given $a>0,b>0$, we denote the constraint by
\beq
T(a,b):=\{(u,v)\in H^1(\R^N)\times H^1(\R^N):\|u\|_2=a,\|v\|_2=b\}.
\eeq
Similar to the scalar case, one can also divide it into several cases, according to whether the energy functional $I$ is bounded from below or not, constrained on $T(a,b)$. For the repulsive case (i.e., $\nu<0$), we refer to Bartsch and Soave\cite{Bartsch2019}. They studied the case of $N=3$ with $p=q=4, \alpha=\beta=2$ and $\mu_1=\mu_2$. For any given $a,b>0$, they could establish the existence of infinitely many positive  normalized solutions. They also studied the phase separation phenomenon as $\nu\rightarrow -\infty$. The study of the attractive case (i.e.,$\nu>0$) is much more tough than the repulsive case.
Based on different situations, let us briefly review the relevant research progress on the attractive case in the following.
\begin{itemize}
\item When $2<p,q,\alpha+\beta<2+\frac{4}{N}$, $I$ is bounded from below on $T(a,b)$ for any $a,b>0$, which is called the mass subcritical case. For this case, Jeanjean and Gou\cite{Gou2016} could obtain the existence of global minimizer, which is also a ground state solution.
\item When $2+\frac{4}{N}<p,q,\alpha+\beta<2^*$, $I$ is unbounded from below on $T(a,b)$ for any $a,b>0$, it is called the mass supercritical case. In particular, it is of pure mass supercritical case. For $N=3$ with $\alpha=\beta=2, p=q=4$, Bartsch, Jeanjean and Soave \cite[Theorem 1.1 and Theorem 1.2]{Bartsch2016} proved the existence of positive normalized solution for $0<\nu<\nu_1$ and $\nu>\nu_2$, where $\nu_1,\nu_2>0$ are defined implicitly by
     $$\max\left\{\frac{1}{a^2\mu_1^2},\frac{1}{b^2\mu_2^2}\right\}=\frac{1}{a^2(\mu_1+\nu_1)^2}+\frac{1}{b^2(\mu_2+\nu_1)^2}$$
     and
     $$\frac{(a^2+b^2)^3}{(\mu_1a^4+\mu_2b^4+2\nu_2a^2b^2)^2}=\min\left\{\frac{1}{a^2\mu_1^2},\frac{1}{b^2\mu_2^2}\right\}.$$
     Clearly the bounds $\nu_1,\nu_2$ depend on the masses $a,b>0$ and
     \beq\lab{eq:property-nu}
     \nu_1\to0,\ \nu_2\to\infty \quad\text{as $\frac{a}{b}\to0$ or $\frac{a}{b}\to\infty$}.
     \eeq
    In particular there is no value of $\nu$ so that the results from \cite{Bartsch2016} yield a solution for all masses. So whether $\nu_1,\nu_2$ can be improved? Furthermore, it is natural to ask that what is the best range of $(a,b,\nu)$ for the existence of positive normalized solution? And this is left as an open problem by Bartsch, Jeanjean and Soave \cite[Remark 1.3-(i)]{Bartsch2016}.
    Bartsch,Zhong and Zou give an answer to Bartsch-Jeanjean-Soave's open problem in \cite{Bartsch2021} and explored almost the best range for the existence of positive normalized solutions, via a new approach different from the framework of variational on the constrained manifold, which is called a global branch approach in the literature now.
    For the more general pure mass supercritical case that $2\leq N\leq 4$ with $\alpha>1,\beta>1, 2+\frac{4}{N}<p,q,\alpha+\beta<2^*$, it was studied by Bartsch and Jeanjean \cite[Theorem 2.4]{Bartsch2018}. For the technical reason, they also obtain the existence result for $0<\nu<\nu_1$ and $\nu>\nu_2$, here $\nu_1,\nu_2$ are also depend on $a,b$ such that \eqref{eq:property-nu} holds. So the Bartsch-Jeanjean-Soave's type open problem is also left here. Recently, a better range of $(a,b,\nu)$ for the existence of positive normalized solutions is obtained by Jeanjean, Zhang and Zhong in \cite{Jeanjean2023}, through a way of variational on the $L^2$-balls.
\item When $2<p,q<2+\frac{4}{N}<\alpha+\beta<2^*$ or $2<\alpha+\beta<2+\frac{4}{N}<p,q<2^*$, it is also a case of mass supercritical since $I$ is also unbounded from below on $T(a,b)$ for any $a,b>0$. However, it is of the mass mixed case here and the geometric structure of the energy functional is much more complicated. Gou and Jeanjean \cite{Gou2018} obtained two different positive normalized solutions provided $\nu>0$ small enough. When $2<p<2+\frac{4}{N}<q<2^*$ and $2+\frac{4}{N}<\alpha+\beta<2^*$, it was considered by Barstch and Jeanjean in \cite{Bartsch2018}, they established the existence of positive normalized solution for some special range of $a,b$.
\end{itemize}

It is worth emphasizing that in the study on normalized solutions for coupled Schr\"odinger system, the following representative literature provides some basic and useful ideas or techniques. The way for constructing a Palais-Smale-Phozaev (PSP for short) sequence, see \cite{Jeanjean1997}. The Pohozave manifold constraint method, see \cite{Bartsch2017} and \cite{Bartsch2019}. The coupled rearrangement, we refer to \cite{Shibata2017} and \cite[Appendix A]{Ikoma2014}. The deformation on the mass constraint manifold, we refer to \cite{Ikoma2019}. A global branch approach to study the normalized solution problem, see \cite{Bartsch2021} and \cite{Jeanjean2021}.

The literature mentioned above do not involve Sobolev critical exponent. But one fact we must admit is that the research on Schr\"odinger equation and system involving Sobolev critical exponent is also a hot topic. Besides the importance in the applications, not negligible reasons of the mathematicians interest for such problems are their stimulating and challenging difficulties coming from the lack of compactness, due to the presence of the limiting exponents for the Sobolev embedding theorems.

For the scalar equation
\beq\lab{eq:20230927-e1}
\begin{cases}
-\Delta u+\lambda u=f(u)~\hbox{in}~\R^N,\\
\int_{\R^N}|u|^2 \ud x=a^2,
\end{cases}
\eeq
the case of $f(u)=\mu |u|^{q-2}u+|u|^{2^*-2}u$ is considered by Soave \cite{Soave2020a}. For a $L^2$-subcritical, $L^2$-critical or $L^2$-supercritical perturbation $\mu |u|^{q-2}u$ of the Sobolev critical term $|u|^{2^*-2}u$, he proved several existence/non-existence and stability/instability results. In particular, for the mass mixed case, i.e., $2<q<2+\frac{4}{N}$, inspired by the result in \cite{Soave2020}, it is natural to believe that there exists at least two different positive normalized solutions for some suitable given mass. Precisely, one is a local minimizer and a second one is a mountain pass solution. However, to recover the compactness of a PSP sequence for the energy functional at the mountain pass level, one has to give a well estimation on the mountain pass value due to the presence of Sobolev critical term. For the technical reason, Soave could not attack it and left it as an open problem, see \cite[Remark 1.1]{Soave2020a}. The key of solving Soave's open problem is finding a suitable testing function to give a better estimation on the mountain pass level.
This open problem is followed and solved by Jeanjean and Le\cite{Jeanjean2022}, Wei and Wu\cite{Wei2022}, Chen and Tang\cite{Chen2022a}. It is worth mentioning that the situation that all the authors of \cite{Soave2020a,Chen2022a,Jeanjean2022a,Wei2022} considered are of double power type. For the general $f(u)$, which is mass mixed and has a Sobolev critical growth, it is studied recently by the authors in \cite{Vicentiu2023}.

For the coupled system case, Bartsch, Li and Zou\cite{Bartsch2023} considered Eq.\eqref{eq:230905-2} with $\mu_1=\mu_2=1$ and $p=q=2^*$. For the case of $\alpha+\beta<2+\frac{4}{N}$, they obtain the existence of positive normalized ground state provided $\nu>0$ small suitable, by searching for a local minimizer. Similar to the scalar case, it is natural to expect a second one solution which is of mountain pass type. However, for the technical reason, the authors in \cite{Bartsch2023} also can not control the mountain pass value to recover the compactness. Hence, Bartsch {\it et al.} raised an open problem focusing on it, see \cite[Remark 1.3-a)]{Bartsch2023}. This open problem can be viewed as a Soave's type open problem for the system case.

In the present paper, we  aim to solve this open problem and this is the main motivation of this project.
We will consider a more general model comparing to that in \cite{Bartsch2023}, but only are interested in the mass mixed case here. More precisely, we shall study the following system
\beq\label{eq:230927-e1}
\begin{cases}
	-\Delta u+\lambda_1 u=\mu_1|u|^{2^*-2}u+\nu\alpha |u|^{\alpha-2}|v|^{\beta}u, ~\hbox{in}~ \R^N,\\
	-\Delta v+\lambda_2 v=\mu_2|v|^{2^*-2}v+\nu\beta |u|^{\alpha}|v|^{\beta-2}v, ~\hbox{in}~ \R^N.
\end{cases}
\eeq
satisfying the normalized constraints \eqref{eq:230814-2}. We only consider the case that
\beq\label{230815-8}
N\in\{3,4\}, a,b>0, ~\hbox{and}~ \alpha,\beta>1, \alpha+\beta<2+\frac{4}{N}.
\eeq
Solutions of the system \eqref{eq:230927-e1}-\eqref{eq:230814-2} can be found as critical points of the functional $I_{\nu}:H^1(\R^N)\times H^1(\R^N)\mapsto \R$,
\beq\label{230904-1}
I_{\nu}(u,v)=\frac{1}{2}\int_{\R^N}|\nabla u|^2\ud x+\frac{1}{2}\int_{\R^N}|\nabla v|^2\ud x-\frac{\mu_1}{2^*}\int_{\R^N}|u|^{2^*}\ud x+\frac{\mu_2}{2^*}\int_{\R^N}|v|^{2^*}\ud x
-\nu\int_{\R^N}|u|^{\alpha}|v|^{\beta}\ud x
\eeq
constrained on the $L^2$ torus $T(a,b)$. For convenience, we denote $H^1(\R^N)\times H^1(\R^N)$ by $E$. We prefer to work on radial subspace and define $E_{rad}=H^1_{rad}(\R^N)\times H_{rad}^1(\R^N)$. From now on, we rewrite $T(a,b)$ by
\beq
T(a,b):=\{(u,v)\in E_{rad}:\|u\|_2=a,\|v\|_2=b\}.
\eeq

Firstly, we are focused on the existence of local minimizers.
\bt\label{th:230820-1}
Suppose \eqref{230815-8} and let $I_{\nu}$ be defined by \eqref{230904-1}. Then the following hold:
\begin{itemize}
\item[(i)] there exists $\nu_0>0$ such that for any fixed $\nu\in(0,\nu_0)$, $I_{\nu}$ has a local minimizer $(u,v)$ on $T(a,b)$, which is positive, radial symmetric and decreasing. The corresponding Lagrange multipliers are positive: $\lambda_1,\lambda_2>0$.
\item[(ii)] Let $\nu_0$ be smaller if necessary, for any $\nu\in(0,\nu_0)$,  the local minimizer $(u,v)$ is indeed a normalized ground state solution to Eqs.\eqref{eq:230927-e1}-\eqref{eq:230814-2}.
\end{itemize}
\et

Secondly, basing on the existence of local minimizer, since $I_{\nu}$ is unbounded from below on $T(a,b)$, it is nature to expect a second one solution which is of mountain pass type. Our next theorem is the highlight of present paper, which gives an affirmative answer to the open problem raised by Bartsch et al.\cite[Remark 1.3-a)]{Bartsch2023}, or we can call it a Soave's type open problem for system.
\bt\label{th:230821-1}
Under the assumptions of Theorem \ref{th:230820-1}-(ii) and let $(u,v)$ be the normalized ground state solution. Then the system \eqref{eq:230927-e1}-\eqref{eq:230814-2} has a second solution $(u_0,v_0)$ in $E_{rad}$, which is of mountain pass type and satisfies
$$I_{\nu}(u,v)<0<I_{\nu}(u_0,v_0)<I_{\nu}(u,v)+\frac{1}{N}\min\{{\mu_1}^{\frac{2-N}{2}},{\mu_2}^{\frac{2-N}{2}}\}S^{\frac{N}{2}}.$$
Moreover, the corresponding Lagrange multipliers are positive: $\bar\lambda_1, \bar\lambda_2>0$.
\et

After proving the existence of multiple normalized solutions, we are also interested in the asymptotic behavior of the solutions obtained by Theorem \ref{th:230820-1} and Theorem \ref{th:230821-1} as $\nu\to0^+$. Indeed, there has been many results on the asymptotic behavior of solutions, see for example \cite{Soave2020,Soave2020a,Wei2022,Jeanjean2022,Bartsch2023,Vicentiu2023}. In particular, for the system \eqref{eq:230927-e1}-\eqref{eq:230814-2}, when $\mu_1=\mu_2=1$, the asymptotic behavior of the normalized ground state solution has been studied in \cite{Bartsch2023}. The system considered here is more general. Indeed, for the mass mixed case, the  asymptotic behavior of local minimizer (ground state) can be studied in a similar way that in \cite{Bartsch2023}. Since the existence of mountain pass solution is new here, our main effort is to capture the  asymptotic behavior of mountain pass solution as $\nu\rightarrow 0$.

For the convenience of describing the results, similar to \cite{Bartsch2023}, we denote by $G(a,b)$ the set of all normalized ground state solutions to the following system
\beq\label{230913-1}
\begin{cases}
-\Delta u+\lambda_1 u=\alpha|u|^{\alpha-2}|v|^{\beta}u ~\hbox{in}~ \R^N,\\
-\Delta v+\lambda_2 v=\beta|u|^{\alpha}|v|^{\beta-2}v ~\hbox{in}~ \R^N,\\
\int_{\R^N}u^2\ud x=a^2, \int_{\R^N}v^2\ud x=b^2,\lambda_1,\lambda_2>0.
\end{cases}
\eeq
Denote by $S$ the sharp constant in the Sobolev embedding, that is,
\beq\label{230913-3}
S=\inf\limits_{u\in D^{1,2}_{0}(\R^N)\setminus \{0\}}\frac{\|\nabla u\|_2^2}{\|u\|_{2^*}^{2}}.
\eeq

We can establish the following result.
\bt\label{th:230913-1}
Suppose \eqref{230815-8} and $\nu\to0^+$. Let $\{(u_{\nu},v_{\nu})\}$ be a family of normalized ground state solutions given by Theorem \ref{th:230820-1}-(ii), $\{(\bar u_{\nu}, \bar v_{\nu})\}$ be a family of mountain pass solutions given by Theorem \ref{th:230821-1}. Let $U$ be the unique positive radial solution of
\beq\label{eq:230921-1}
-\Delta u=u^{2^*-1} ~\hbox{in}~ \R^N.
\eeq
Then the following hold:
\begin{itemize}
\item [(i)] $\|u_{\nu}\|_{\infty}+\|v_{\nu}\|_{\infty}\to0$ and there exists $t_{\nu}>0$, $t_{\nu}\sim\nu^{\frac{1}{2-\gamma_{\alpha+\beta}}}$ such that
$$dist_{E}\left(t_{\nu}\star (u_{\nu},v_{\nu}),G(a,b)\right)\to0, ~\hbox{as}~ \nu\to 0^+,$$
where $t\star (u,v)=(t\star u,t\star v):=(t^{\frac{N}{2}}u(tx), t^{\frac{N}{2}}v(tx))$.
\item [(ii-1)] If $\mu_1<\mu_2$, then there exists $\varepsilon_{\nu}>0$ such that $(\bar u_{\nu},  \varepsilon_{\nu}^{\frac{N-2}{2}}\bar v_{\nu}(\varepsilon_{\nu}x))\to(0,\mu_{2}^{\frac{N-2}{4}}U)$ in $D_{0}^{1,2}(\R^N)\times D_{0}^{1,2}(\R^N)$ as $\nu\to0^+$.
\item [(ii-2)] If $\mu_2<\mu_1$, then there exists $\varepsilon_{\nu}>0$ such that $(\varepsilon_{\nu}^{\frac{N-2}{2}}\bar u_{\nu}(\varepsilon_{\nu}x), \bar v_{\nu})\to(\mu_{1}^{\frac{N-2}{4}}U,0)$ in $D_{0}^{1,2}(\R^N)\times D_{0}^{1,2}(\R^N)$ as $\nu\to0^+$.
\item [(ii-3)] If $\mu_1=\mu_2$, then one of the statements in (ii-1) or (ii-2) holds.
\end{itemize}
\et

The paper is organized as follows. In section \ref{sec:Preliminaries}, we prove the existence of the local minimizer and Theorem \ref{th:230820-1}. In section \ref{sec:3}, we mainly prove Theorem \ref{th:230821-1}. In section \ref{sec:4}, we focus on the asymptotic behavior of the solutions and complete the proof of Theorem \ref{th:230913-1}.

\s{The existence of local minimizer and Proof of Theorem \ref{th:230820-1}}\lab{sec:Preliminaries}
\subsection{Local minimal structure }
\bl\label{le:230903-1}
Suppose \eqref{230815-8}. There exists $\nu_0, k_0, \rho_0>0$ such that for any $\nu<\nu_0$,
\beq\label{230902-2}
m_{\nu}(a,b):=\inf_{V_{\rho_0}}I_{\nu}(u,v)<0<k_0\leq \inf_{\partial V_{\rho_0}}I_{\nu}(u,v),
\eeq
where $V_{\rho_0}$ and $\partial V_{\rho_0}$ are given by
$$V_{\rho_0}:=\{(u,v)\in T(a,b):(\|\nabla u\|_2^2+\|\nabla v\|_2^2)^{\frac{1}{2}}<\rho_0\}$$
and
$$\partial V_{\rho_0}:=\{(u,v)\in T(a,b):(\|\nabla u\|_2^2+\|\nabla v\|_2^2)^{\frac{1}{2}}=\rho_0\}.$$
\el
\bp
For any $(u,v)\in\partial V_{\rho_0}$, where $\rho_0>0$ to be determined, by Gagliardo-Nirenberg inequality, i.e., there exists $C(N,\alpha,\beta)>0$ such that
\beq\label{230901-1}
\||u|^{\alpha}|v|^{\beta}\|_1\leq C(N,\alpha,\beta)(\|u\|_2^2+\|v\|_2^2)^{(\alpha+\beta-\gamma_{\alpha+\beta})/2}(\|\nabla u\|_2^2+\|\nabla v\|_2^2)^{\gamma_{\alpha+\beta}/2},
\eeq
combining with critical Sobolev inequality,
it follows that
\beq\label{230902-1}
\begin{aligned}
I_{\nu}(u,v)=&\frac{1}{2}(\|\nabla u\|_2^2+\|\nabla v\|_2^2)-\frac{\mu_1}{2^*}\|u\|_{2^*}^{2^*}-\frac{\mu_2}{2^*}\|v\|_{2^*}^{2^*}-\nu\int_{\R^N}|u|^{\alpha}|v|^{\beta}\ud x\\
\geq&\frac{1}{2}(\|\nabla u\|_2^2+\|\nabla v\|_2^2)-\frac{\mu_1}{2^*}S^{-\frac{2^*}{2}}\|\nabla u\|_{2}^{2^*}-\frac{\mu_2}{2^*}S^{-\frac{2^*}{2}}\|\nabla v\|_{2}^{2^*}\\
&-\nu C(N,\alpha,\beta)(a^2+b^2)^{(\alpha+\beta-\gamma_{\alpha+\beta})/2}(\|\nabla u\|_2^2+\|\nabla v\|_2^2)^{\gamma_{\alpha+\beta}/2}\\
\geq&\frac{1}{2}\rho_0^2-\frac{1}{2^*}S^{-\frac{2^*}{2}}\max\{\mu_1,\mu_2\}\rho_{0}^{2^*}-\nu C(N,\alpha,\beta)(a^2+b^2)^{(\alpha+\beta-\gamma_{\alpha+\beta})/2}\rho_0^{\gamma_{\alpha+\beta}}\\
=&\rho_0^2(\frac{1}{2}-\frac{1}{2^*}S^{-\frac{2^*}{2}}\max\{\mu_1,\mu_2\}\rho_{0}^{2^*-2}-\nu C(N,\alpha,\beta)(a^2+b^2)^{(\alpha+\beta-\gamma_{\alpha+\beta})/2}\rho_0^{\gamma_{\alpha+\beta}-2}).
\end{aligned}
\eeq
For $\rho>0$, we define
$$
h_{\nu}(\rho):=\frac{1}{2}-\nu A\rho^{\gamma_{\alpha+\beta}-2}-B\rho^{2^*-2},
$$
where
$$A:= C(N,\alpha,\beta)(a^2+b^2)^{\frac{\alpha+\beta-\gamma_{\alpha+\beta}}{2}}~\hbox{and}~ B:=\frac{1}{2^*}S^{-\frac{2^*}{2}}\max\{\mu_1,\mu_2\}.$$
One can obtain that $h_{\nu}'(\rho)=-\nu A(\gamma_{\alpha+\beta}-2)\rho^{\gamma_{\alpha+\beta}-3}-B(2^*-2)\rho^{2^*-3}$ and by $h_{\nu}'(\rho)=0$ there exists a unique solution:
\beq\lab{eq:20230908-e1}
\rho_{\nu}=\left[\frac{\nu A(2-\gamma_{\alpha+\beta})}{B(2^*-2)}\right]^{\frac{1}{2^*-\gamma_{\alpha+\beta}}}=:\nu^{{\frac{1}{2^*-\gamma_{\alpha+\beta}}}}C_1.
\eeq

Noting that $h_{\nu}(\rho)\to-\infty$ provided $\rho\to 0$ or $\rho\to+\infty$, so $$
\begin{aligned}
\max\limits_{\rho>0}h_{\nu}(\rho)&=h_{\nu}(\rho_{\nu})=\frac{1}{2}-(AC_{1}^{\gamma_{\alpha+\beta}-2}+BC_{1}^{2^*-2})\nu^{\frac{2^*-2}{2^*-\gamma_{\alpha+\beta}}}\\
&=\frac{1}{2}-\left((\frac{2-\gamma_{\alpha+\beta}}{2^*-2})^{\frac{\gamma_{\alpha+\beta}-2}{2^*-\gamma_{\alpha+\beta}}}+(\frac{2-\gamma_{\alpha+\beta}}{2^*-2})^{\frac{2^*-2}{2^*-\gamma_{\alpha+\beta}}}\right)A^{\frac{2^*-2}{2^*-\gamma_{\alpha+\beta}}}B^{\frac{2-\gamma_{\alpha+\beta}}{2^*-\gamma_{\alpha+\beta}}}\nu^{\frac{2^*-2}{2^*-\gamma_{\alpha+\beta}}}.
\end{aligned}
$$
By $\max\limits_{\rho>0}h_{\nu}(\rho)>0$ we have that $\nu<\left(\frac{1}{2(AC_{1}^{\gamma_{\alpha+\beta}-2}+BC_{1}^{2^*-2})}\right)^{\frac{2^*-\gamma_{\alpha+\beta}}{2^*-2}}=:\bar\nu_0$, i.e.
\beq\label{230902-3}
\bar\nu_0=(\max\{\mu_1,\mu_2\})^{\frac{\gamma_{\alpha+\beta}-2}{2^*-2}}\frac{(2^*-2)(2^*(2-\gamma_{\alpha+\beta}))^{\frac{2-\gamma_{\alpha+\beta}}{2^*-2}}}{(2(2^*-\gamma_{\alpha+\beta}))^{\frac{2^*-\gamma_{\alpha+\beta}}{2^*-2}}}\cdot \frac{S^{\frac{2^*(2-\gamma_{\alpha+\beta})}{2(2^*-2)}}}{C(N,\alpha,\beta)(a^2+b^2)^{\frac{\alpha+\beta-\gamma_{\alpha+\beta}}{2}}}.
\eeq
Now we let $\rho_0:=\bar\nu_{0}^{{\frac{1}{2^*-\gamma_{\alpha+\beta}}}}C_1$, then for any $\nu<\bar{\nu}_0$, we have $h_{\nu}(\rho_0)>h_{\bar\nu_0}(\rho_0)=0$. We choose $\nu_0<\bar\nu_0$ and let $k_0:=\rho_0^2h_{\nu_0}(\rho_0)>0$. Hence for any $\nu\in(0,\nu_0)$, by \eqref{230902-1} we obtain $I_{\nu}(u,v)\geq \rho_0^2 h_{\nu}(\rho_0)>k_0>0$. By the arbitrariness of $(u,v)\in\partial V_{\rho_0}$, we have $\inf\limits_{\partial V_{\rho_0}}I_{\nu}(u,v)\geq k_0>0$.

Take $(u,v)\in T(a,b)$ such that $\int_{\R^N}|u|^{\alpha}|v|^{\beta}\ud x\neq0$ and let
\beq
t\star (u,v)=(t\star u,t\star v):=(t^{\frac{N}{2}}u(tx), t^{\frac{N}{2}}v(tx)), t>0,
\eeq
then $(t\star u, t\star v)\in V_{\rho_0}$ when $t$ is small enough. And by a computation, we have
\beq\label{230903-1}
I_{\nu}(t\star u,t\star v)=\frac{t^2}{2}(\|\nabla u\|_2^2+\|\nabla v\|_2^2)-\frac{\mu_1 t^{2^*}}{2^*}\|u\|_{2^*}^{2^*}-\frac{\mu_2 t^{2^*}}{2^*}\|v\|_{2^*}^{2^*}-\nu t^{\frac{(\alpha+\beta-2)N}{2}}\int_{\R^N}|u|^{\alpha}|v|^{\beta}\ud x,
\eeq
which implies that $I_{\nu}(t\star u,t\star v)<0$ provided $t$ is small enough since $\frac{(\alpha+\beta-2)N}{2}<2$. Hence we obtain \eqref{230902-2}.
\ep
\br\label{remark:230902-1}
 The formula \eqref{230902-3} declares that the sharp value of $\bar\nu_0$ depends on the masses $a$ and $b$. Let us denote it by $\bar\nu_0(a,b)$. Then it is easy to check that $\bar\nu_0(a,b)$ decreases in the sense that $\bar\nu_0(a_1,b_1)\geq \bar\nu_0(a,b)$ provided $a_1\leq a$ and $b_1\leq b$. Hence, for any $a_1\leq a$,$b_1\leq b$ with $\nu<\bar\nu_0(a,b)$, it is true that $\nu<\bar\nu_0(a_1,b_1)$. And thus \eqref{230902-2} keeps valid after replacing $m_{\nu}(a,b)$ by $m_{\nu}(a_1,b_1)$ and $V_{\rho_0}$  is modified by
$V_{\rho_0}=\{(u,v)\in T(a_1,b_1):(\|\nabla u\|_2^2+\|\nabla v\|_2^2)^{\frac{1}{2}}<\rho_0\}$.
\er

\br\label{br:230913-1}
For $0<\nu<\nu_0$, suppose that $(u_\nu,v_\nu)$ is a local minimizer, then by the definition of $\rho_\nu$, it is trivial that $(\|\nabla u_\nu\|_2^2+\|\nabla v_\nu\|_2^2)^{\frac{1}{2}}<\rho_\nu$. So by \eqref{eq:20230908-e1}, one can see that $\rho_\nu\rightarrow 0$ as  $\nu\rightarrow 0$ and thus
$(u_\nu,v_\nu)\rightarrow (0,0)$ in $D_{0}^{1,2}(\R^N)\times D_{0}^{1,2}(\R^N)$.
\er

\subsection{Existence of local minimizer}
To prove the existence of local minimizer, we introduce a new set:
$$\widetilde V_{\rho_0}:=\{(u,v)\in E_{rad}:\|u\|_2\leq a, \|v\|_2\leq b,\|\nabla u\|_2^2+\|\nabla v\|_2^2\leq \rho_0^2\}.$$
Setting $\tilde m_{\nu}(a,b):=\inf\limits_{\widetilde V_{\rho_0}}I_{\nu}(u,v)$, then it is trivial that $\tilde m_{\nu}(a,b)\leq m_{\nu}(a,b)<0$. Firstly, we prove that $\tilde m_{\nu}(a,b)$ is always achievable.

\bl\label{le:230902-1}
Suppose \eqref{230815-8}  and let $\nu_0>0,\rho_0>0$ be the numbers given by Lemma \ref{le:230903-1}.
 Then for any $\nu\in (0,\nu_0)$, $\tilde m_{\nu}(a,b)$ is attained by some $(u,v)\in\widetilde V_{\rho_0}$ and $u\gneqq0$, $v\gneqq0$.
\el
\bp
Take a minimizing sequence $\{(u_n,v_n)\}\subset\widetilde V_{\rho_0}$. Denote the Schwartz rearrangement of $u$ by $u^*$, then it is well known that $\|\nabla u^*\|_2^2\leq\|\nabla u\|_2^2$, $\|u^*\|_2^2=\|u\|_2^2$, $\|u^*\|_{2^*}^{2^*}=\|u\|_{2^*}^{2^*}$ and $\int_{\R^N}|u^*|^{\alpha}|v^*|^{\beta}\ud x\geq\int_{\R^N}|u|^{\alpha}|v|^{\beta}\ud x$. Hence $\{(u_{n}^{*},v_n^*)\}$ is also a minimizing sequence for $I_{\nu}$ on $\widetilde V_{\rho_0}$. Without loss of generality, we may assume that $(u_n,v_n)=(u_n^*,v_n^*)$, then $\{u_n\}$ and $\{v_n\}$ are nonnegative radial non-increasing functions.

Since $\{(u_n,v_n)\}\subset\widetilde V_{\rho_0}$ and $\widetilde V_{\rho_0}$ is bounded, closed and convex, we conclude that $\{(u_n,v_n)\}$ is bounded in $E_{rad}$ and we may assume that $(u_n,v_n)\rightharpoonup (u,v)$ in $E_{rad}$, then $(u,v)\in\widetilde V_{\rho_0}$. We {\bf claim} that $u\not\equiv0$ and $v\not\equiv0$. Suppose $u\equiv0$ or $v\equiv0$, then $\int_{\R^N}|u_n|^{\alpha}|v_n|^{\beta}\ud x\leq\|u_n\|_{\alpha+\beta}^{\alpha}\|v_n\|_{\alpha+\beta}^{\beta}\to0$. Up to a subsequence, we assume that $\displaystyle\lim_{n\rightarrow \infty}(\|\nabla u_n\|_2^2+\|\nabla v_n\|_2^2)=\rho^2$. Then we have $\rho^2\leq\rho_0^2$ and
\begin{align*}
\tilde{m}_\nu(a,b)=&I_{\nu}(u_n,v_n)+o_n(1)\\
=&\frac{1}{2}(\|\nabla u_n\|_2^2+\|\nabla v_n\|_2^2)-\frac{\mu_1}{2^*}\|u_n\|_{2^*}^{2^*}-\frac{\mu_2}{2^*}\|v_n\|_{2^*}^{2^*}-\nu\int_{\R^N}|u_n|^{\alpha}|v_n|^{\beta}\ud x+o_n(1)\\
\geq&\frac{1}{2}(\|\nabla u_n\|_2^2+\|\nabla v_n\|_2^2)-\frac{\mu_1}{2^*}S^{-\frac{2^*}{2}}\|\nabla u_n\|_{2}^{2^*}-\frac{\mu_2}{2^*}S^{-\frac{2^*}{2}}\|\nabla v_n\|_{2}^{2^*}+o_n(1)\\
\geq&\frac{1}{2}\rho^2-\frac{1}{2^*}S^{-\frac{2^*}{2}}\max\{\mu_1,\mu_2\}\rho^{2^*}+o_n(1)\\
\geq&\rho^2(\frac{1}{2}-\frac{1}{2^*}S^{-\frac{2^*}{2}}\max\{\mu_1,\mu_2\}\rho_{0}^{2^*-2})+o_n(1)\geq0,
\end{align*}
where the last inequality follows from the fact of $h_{\nu}(\rho_0)>h_{\nu_0}(\rho_0)>0$. This contradicts with $\tilde m_{\nu}(a,b)<0$ and the claim is proved.

Since $\{(u_n,v_n)\}$ is bounded in $E_{rad}$, by Br\'ezis-Lieb Lemma, for $n$ large enough, we have
$$\|\nabla (u_n-u)\|_2^2+\|\nabla (v_n-v)\|_2^2=\|\nabla u_n\|_2^2+\|\nabla v_n\|_2^2-\|\nabla u\|_2^2-\|\nabla v\|_2^2+o_n(1)<\rho_0^2,$$
$$\|u_n-u\|_2^2=\|u_n\|_2^2-\|u\|_2^2+o_n(1)< a^2 ~\hbox{and}~ \|v_n-v\|_2^2=\|v_n\|_2^2-\|v\|_2^2+o_n(1)< b^2,$$
which implies that $(u_n-u,v_n-v)\in \widetilde V_{\rho_0}$. By Br\'ezis-Lieb Lemma again, we have
\begin{align*}
&\tilde{m}_{\nu}(a,b)+o_n(1)=I_{\nu}(u_n,v_n)\\
=&I_{\nu}(u,v)+\frac{1}{2}(\|\nabla (u_n-u)\|_2^2+\|\nabla (v_n-v)\|_2^2)-\frac{\mu_1}{2^*}\|u_n-u\|_{2^*}^{2^*}-\frac{\mu_2}{2^*}\|v_n-v\|_{2^*}^{2^*}+o_n(1)\\
\geq&\tilde{m}_{\nu}(a,b)+\frac{1}{2}(\|\nabla (u_n-u)\|_2^2+\|\nabla (v_n-v)\|_2^2)-\frac{\mu_1}{2^*}S^{-\frac{2^*}{2}}\|\nabla (u_n-u)\|_{2}^{2^*}\\
&-\frac{\mu_2}{2^*}S^{-\frac{2^*}{2}}\|\nabla (v_n-v)\|_{2}^{2^*}+o_n(1)\\
\geq&\tilde{m}_{\nu}(a,b)+o_n(1),
\end{align*}
which implies that $I_{\nu}(u,v)=\tilde m_{\nu}(a,b)$. Hence $(u,v)$ obtains $\tilde m_{\nu}(a,b)$ and $u\gneqq0$, $v\gneqq0$.
\ep

\br\label{re:230902-2}
Under the assumptions of Lemma \ref{le:230902-1}, for any local minimizer $(u,v)$, i.e., $(u,v)\in \tilde{V}_{\rho_0}$ and $I_\nu(u,v)=\tilde m_{\nu}(a,b)$, we claim that $\|\nabla u\|_2^2+\|\nabla v\|_2^2<\rho_0^2$. Indeed, if $\|\nabla u\|_2^2+\|\nabla v\|_2^2=\rho_0^2$, then by Remark \ref{remark:230902-1} we have $I_{\nu}(u,v)>0$, which contradicts to the fact that $\tilde m_{\nu}(a,b)<0$. Hence,
$$
\tilde m_{\nu}(a,b)=\min\{I_{\nu}(u,v):(u,v)\in E_{rad},\|u\|_2\leq a, \|v\|_2\leq b,\|\nabla u\|_2^2+\|\nabla v\|_2^2<\rho_0^2\}.
$$
Thus, there exist some $\lambda_1,\lambda_2\in\R$ such that $(u,v)$ is a couple nontrivial nonnegative radial decreasing solution to
\beq\label{230902-4}
\begin{cases}
-\Delta u+\lambda_1 u=\mu_1 u^{2^*-1}+\nu\alpha u^{\alpha-1}v^{\beta} ~\hbox{in}~ \R^N,\\
-\Delta v+\lambda_2 v=\mu_2 v^{2^*-1}+\nu\beta u^{\alpha}v^{\beta-1} ~\hbox{in}~ \R^N.\\
\end{cases}
\eeq
In particular, by the standard elliptic regularity, one can prove that $u,v$ are smooth.
\er
To make clear that whether $(u,v)\in T(a,b)$ or not, inspired by the scalar case in \cite[Lemma 4.5]{Vicentiu2023}, we can establish the following result:
\bl\lab{lemma:20230907-l1}
Let $(u,v)$ and $\lambda_1,\lambda_2$ be given in Remark \ref{re:230902-2}. The following are true:
\begin{itemize}
\item [(i)] If $\|u\|_2<a$, then $\lambda_1=0$; If $\|v\|_2<b$, then $\lambda_2=0$.
\item [(ii)] If $\|u\|_2=a$, then $\lambda_1\geq0$; If $\|v\|_2=b$, then $\lambda_2\geq0$.
\end{itemize}
\el
\bp
The argument is the same as that in \cite[Lemma 4.5]{Vicentiu2023}. One only needs to replace $I(u)$ by $I_{\nu}(u,v)$. We omit the details here.
\ep

\begin{corollary}\label{corollary:230902-1}
Let $(u,v)$ and $\lambda_1,\lambda_2$ be given in Remark \ref{re:230902-2}.
 If $\lambda_1\neq0$ and $\lambda_2\neq0$, then $(u,v)\in T(a,b)$ and $\lambda_1>0$, $\lambda_2>0$.
\end{corollary}
\bp
This is just a consequence follows directly from Lemma \ref{lemma:20230907-l1}.
\ep

\bl\label{le:230821-1}
Under the assumptions of Lemma \ref{le:230902-1},
let $\tilde m_{\nu}(a,b),m_{\nu}(a,b)$ be defined above. Then for any $0<\nu<\nu_0$, the following statements are true.
\begin{itemize}
\item [(i)] $\tilde m_{\nu}(a,b)=m_{\nu}(a,b)$ and  $m_{\nu}(a,b)$ is achievable.
\item [(ii)] $m_{\nu}(a,b)\leq m_{\nu}(a_1,b_1)$ for any $0\leq a_1\leq a, 0\leq b_1\leq b$. Moreover, the equality holds if and only if $(a_1,b_1)=(a,b)$.
\end{itemize}
\el
\bp
$(i)$ By Lemma \ref{le:230902-1}, there exist some $(u,v)\in\widetilde V_{\rho_0}$ such that $I_{\nu}(u,v)=\tilde m_{\nu}(a,b)$. If $(u,v)\in T(a,b)$, by Remark \ref{re:230902-2} we know that $(u,v)\in V_{\rho_0}$ and then $m_{\nu}(u,v)\leq I_{\nu}(u,v)=\tilde m_{\nu}(a,b)\leq m_{\nu}(a,b)$, which implies $\tilde m_{\nu}(a,b)=m_{\nu}(a,b)$ and thus $(u,v)$ is a local minimizer of $I_{\nu}$ on $T(a,b)$. To prove $(u,v)\in T(a,b)$, by Corollary \ref{corollary:230902-1}, we just need to show $\lambda_1\neq0$ and $\lambda_2\neq0$. Suppose $\lambda_1=0$, then by \eqref{230902-4} we have that
$$-\Delta u=\mu_1 u^{2^*-1}+\nu\alpha u^{\alpha-1}v^{\beta}\geq0 ~\hbox{in}~ \R^N.$$
By Remark \ref{re:230902-2}, $u$ is smooth. Then by Liouville theorem \cite[Lemma A.2]{Ikoma2014}, we conclude that $u\equiv0$, which is a contradiction. Similarly, we can show that $\lambda_2\neq0$.

$(ii)$ For  any $0<a_1\leq a, 0<b_1\leq b$, recalling Remark \ref{remark:230902-1}, the conclusion of (i) declares that $m_{\nu}(a_1,b_1)=\tilde m_{\nu}(a_1,b_1)$ and it is achievable. So,
 $m_{\nu}(a_1,b_1)=\tilde m_{\nu}(a_1,b_1)\geq\tilde m_{\nu}(a,b)= m_{\nu}(a,b)$.
On the other hand, by the choice of $\rho_0$, for any $(u,v)\in \tilde{V}_{\rho_0}$ and at least one of $u,v$ is trivial, one can check that $I_\nu(u,v)\geq 0$. So, if $a_1=0$ or $b_1=0$, it holds that $m_\nu(a_1,b_1)\geq 0$. Hence, we also have
$m_\nu(a,b)<0\leq m_\nu(a_1,b_1)$. Suppose $(a_1,b_1)\neq (a,b)$, we may assume that $a_1<a, b_1\leq b$. If $m_{\nu}(a_1,b_1)=m_{\nu}(a,b_1)$ and assume that $(u,v)\in T(a_1,b_1)$ attains $m_{\nu}(a_1,b_1)$, then by Remark \ref{re:230902-2} we know that $(u,v)$ is a solution to \eqref{230902-4} and by Lemma \ref{le:230821-1} (i) we have $\lambda_1\neq0$. On the other hand, by Lemma \ref{lemma:20230907-l1} (i), we conclude that $\lambda_1=0$, which is a contradiction. Thus $m_{\nu}(a_1,b_1)>m_{\nu}(a,b_1)\geq m_{\nu}(a,b)$.
\ep
\vskip 0.2in
\noindent
{\bf Proof of Theorem \ref{th:230820-1}:}
Combining with Lemma \ref{le:230902-1} and Lemma \ref{le:230821-1} we have that $(u,v)$ is a local minimizer of $I_{\nu}$ on $T(a,b)$ and $\lambda_1>0,\lambda_2>0$.
Next we prove that there exists $\nu_0$ small such that for any $\nu\in(0,\nu_0)$, the local minimizer $(u,v)$ is also a normalized ground state. If the statement is false, then there exists a sequence $\nu_n\to0$, assume that $(u_n,v_n)$ is the local minimizer of $I_{\nu_n}$ on $T(a,b)$ and there exists another critical point $(\tilde u_n,\tilde v_n)$ of $I_{\nu_n}$ on $T(a,b)$ such that $I_{\nu_n}(\tilde u_n,\tilde v_n)<I_{\nu_n}(u_n,v_n)=m_{\nu_n}(a,b)<0$.
Note that a solution must satisfy the so-called pohozaev identity $P_{\nu_n}(\tilde u_n,\tilde v_n)=0$, where the definition of $P_\nu(u,v)$ we refer to \eqref{230904-2}.
Hence, combining with \eqref{230901-1}, we have
\begin{align*}
0>&I_{\nu_n}(\tilde u_n,\tilde v_n)=\frac{1}{2}(\|\nabla\tilde u_n\|_2^2+\|\nabla\tilde v_n\|_2^2)-\frac{\mu_1}{2^*}\|\tilde u_n\|_{2^*}^{2^*}-\frac{\mu_2}{2^*}\|\tilde v_n\|_{2^*}^{2^*}-\nu_n\int_{\R^N}|\tilde u_n|^{\alpha}|\tilde v_n|^{\beta}\ud x\\
=&(\frac{1}{2}-\frac{1}{2^*})(\|\nabla\tilde u_n\|_2^2+\|\nabla\tilde v_n\|_2^2)-\nu_n\frac{2^*-\gamma_{\alpha+\beta}}{2^*}\int_{\R^N}|\tilde u_n|^{\alpha}|\tilde v_n|^{\beta}\ud x\\
\geq&(\frac{1}{2}-\frac{1}{2^*})(\|\nabla\tilde u_n\|_2^2+\|\nabla\tilde v_n\|_2^2)-\nu_n\frac{2^*-\gamma_{\alpha+\beta}}{2^*}(a^2+b^2)^{\frac{\alpha+\beta-\gamma_{\alpha+\beta}}{2}}(\|\nabla\tilde u_n\|_2^2+\|\nabla\tilde v_n\|_2^2)^{\frac{\gamma_{\alpha+\beta}}{2}},
\end{align*}
which implies that
$$0\leq(\frac{1}{2}-\frac{1}{2^*})(\|\nabla\tilde u_n\|_2^2+\|\nabla\tilde v_n\|_2^2)^{\frac{2-\gamma_{\alpha+\beta}}{2}}<\nu_n\frac{2^*-\gamma_{\alpha+\beta}}{2^*}(a^2+b^2)^{\frac{\alpha+\beta-\gamma_{\alpha+\beta}}{2}}.$$
So $\|\nabla\tilde u_n\|_2^2+\|\nabla\tilde v_n\|_2^2\to0$ as $\nu_n\to0$ and thus $(\tilde u_n,\tilde v_n)\in V_{\rho_0}$ for $n$ large enough. Hence, we have $I_{\nu_n}(\tilde u_n,\tilde v_n)\geq m_{\nu_n}(a,b)$, this is a contradiction.
\hfill$\Box$

\s{the existence of the mountain pass solution}\label{sec:3}
For the Sobolev critical case, $I_{\nu}$ is unbounded below from $T(a,b)$, basing on the existence of local minimizer $(u,v)$, it is expected that $I_\nu|_{T(a,b)}$ has a mountain pass geometric near $(u,v)$. In this section, we devote to solve the open problem raised by Bartsch et al.\cite[Remark 1.3-a)]{Bartsch2023}. Precisely, under the assumptions of Theorem \ref{th:230820-1}-(ii), we shall prove the existence of mountain pass solution, which is exactly the main content of the open problem.
\subsection{Mountain pass geometric structure}
\bl\label{le:230820-1}
Under the assumptions of Theorem \ref{th:230820-1}, let $(u,v)$ be the local minimizer and $k_0$ be given by Lemma \ref{le:230903-1}, then
\beq\label{230921-1}
M_{\nu}(a,b):=\inf\limits_{\gamma\in\Gamma_{\nu}}\max\limits_{t\in[0,1]}I_{\nu}(\gamma(t))\geq k_0>\sup\limits_{\gamma\in\Gamma_{\nu}}\max\{I_{\nu}(\gamma(0)),I_{\nu}(\gamma(1))\},
\eeq
where
$$\Gamma_{\nu}=\{\gamma\in C([0,1],T(a,b)):\gamma(t)=(\gamma_1(t),\gamma_2(t)),\gamma(0)=(u,v),I_{\nu}(\gamma(1))<2m_{\nu}(a,b)\}.$$
\el
\bp
Firstly, we remark that $\Gamma_{\nu}\neq\emptyset$ and thus $M_{\nu}(a,b)$ is well-defined. Indeed, recalling the fiber map $t\star (u,v)=(t^{\frac{N}{2}}u(tx),t^{\frac{N}{2}}v(tx))$ and \eqref{230903-1}, we know that $t\star (u,v)\in T(a,b)$ and there exists $T>0$ large such that $I_{\nu}(T\star u,T\star v)<2m_{\nu}(a,b)<0$. Letting $\gamma(t):=(1+(T-1)t)\star (u,v)$, then $\gamma(t)\in\Gamma_{\nu}$.

By Lemma \ref{le:230903-1}, we know that $k_0=\rho_0^2 h_{\nu_0}(\rho_0)>0>\sup\limits_{\gamma\in\Gamma_{\nu}}\max\{I_{\nu}(\gamma(0)),I_{\nu}(\gamma(1))\}$. For any $\gamma\in\Gamma_{\nu}$, since $I_{\nu}(\gamma(1))<2m_{\nu}(a,b)$, we know that $\gamma(1)\not\in \overline V_{\rho_0}$, i.e. $(\|\nabla \gamma_1(1)\|_2^2+\|\nabla\gamma_2(1)\|_2^2)^{\frac{1}{2}}>\rho_0$. Noting that $(\|\nabla \gamma_1(0)\|_2^2+\|\nabla\gamma_2(0)\|_2^2)^{\frac{1}{2}}<\rho_0$, by Intermediate Value Theorem, there exists $t_0\in(0,1)$ such that
$(\|\nabla\gamma_1(t_0)\|_2^2+\|\nabla\gamma_2(t_0)\|_2^2)^{\frac{1}{2}}=\rho_0$, which implies that $\gamma(t_0)\in\partial V_{\rho_0}$. Hence, by Lemma \ref{le:230903-1}, we have $\max\limits_{t\in[0,1]}I_{\nu}(\gamma(t))\geq I_{\nu}(\gamma(t_0))\geq\inf\limits_{\partial V_{\rho_0}}I_{\nu}(u,v)\geq k_0$. By the arbitrariness of $\gamma$, we conclude the desired result.
\ep

\subsection{Existence of nonnegative $(PSP)_{M_{\nu}(a,b)}$ sequence}
Since the situation we considered is of mass supercritical case, a Palais-Smale sequence is not necessary bounded in $E$.
In order to find a bounded Palais-Smale sequence, we expect to find a (PSP) sequence, which is a Palais-Smale sequence with an additional property $P_\nu(u_n,v_n)\rightarrow 0$. Here $P_\nu(u,v)$ is defined by
\beq\label{230904-2}
P_{\nu}(u,v)=\|\nabla u\|_2^2+\|\nabla v\|_2^2-\mu_1\|u\|_{2^*}^{2^*}-\mu_2\|v\|_{2^*}^{2^*}-\nu \gamma_{\alpha+\beta}\int_{\R^N}|u|^{\alpha}|v|^{\beta}dx,
\eeq
where $\gamma_{\alpha+\beta}=\frac{N(\alpha+\beta-2)}{2}$. For this goal, we shall adopt the idea by introducing an auxiliary functional $\tilde I_{\nu}$ in \cite{Jeanjean1997} and apply the mini-max principle in \cite{Ghoussoub1993} (see also Lemma \ref{le:230820-2}), one can also refer to \cite{Ikoma2019}.

Recalling that $E=H^1(\R^N)\times H^1(\R^N)$ and letting $H=L^2(\R^N)\times L^2(\R^N)$, for any $(u,v)\in E$, we define the norms of $E$ and $H$ by
$$
\|(u,v)\|^{2}_{E}:=\|u\|_{H^1(\R^N)}^{2}+\|v\|_{H^1(\R^N)}^{2}, \quad \|(u,v)\|^{2}_{H}:=\|u\|_2^2+\|v\|_2^2.
$$
Define a continuous map $\Theta:H^1(\R^N)\times \R\to H^1(\R^N)$ by
$$
\Theta(v,s)(x):=e^{\frac{Ns}{2}}v(e^sx) ~\hbox{for}~ v\in E, \forall s\in \R, x\in{\R^N}.
$$
We introduce the auxiliary functional $\tilde I_{\nu}$:
\beq\label{230820-5}
\begin{aligned}
\tilde I_{\nu}(s,u,v):&=I_{\nu}(\Theta(u,s),\Theta(v,s))=\frac{e^{2s}}{2}(\|\nabla u\|_2^2+\|\nabla v\|_2^2)-\frac{e^{2^*s}}{2^*}(\mu_1\|u\|_{2^*}^{2^*}+\mu_2\|v\|_{2^*}^{2^*})\\
&-\nu e^{(\frac{\alpha+\beta}{2}-1)Ns}\int_{\R^N}|u|^{\alpha}|v|^{\beta}\ud x.
\end{aligned}
\eeq
By some direct computations, similar to \cite{Ikoma2019,Chen2022a}, for any $(s,u,v)\in \R\times E$, one can deduce the following two formulas,
\beq\label{230820-3}
P_{\nu}(\Theta(u,s),\Theta(v,s))=|\langle \tilde I_{\nu}'(s,u,v),(1,0,0)\rangle|\leq\|\tilde I_{\nu}|'_{\R\times T(a,b)}(s,u,v)\|,
\eeq
\beq\label{230820-4}
\|I_{\nu}|'_{T(a,b)}(\Theta(u,s),\Theta(v,s))\|
\leq e^{|s|}\|\tilde I_{\nu}|'_{\R\times T(a,b)}(s,u,v)\|.
\eeq

To obtain a $(PSP)$ sequence $\{(u_n,v_n)\}$, exactly a $(PS)$ sequence with an additional condition $P_{\nu}(u_n,v_n)\to0$, the following lemma is usually essential for us.
\bl\label{le:230820-2}
(cf.\cite[Theorem 3.2]{Ghoussoub1993})
Let $\varphi$ be a $C^1$-functional on a complete connected $C^1$-Finsler manifold $X$ (without boundary) and consider a homotopy-stable family $\cf$ of compact subsets of $X$ with a closed boundary $B$. Set $c=c(\varphi,\cf)=\inf\limits_{A\in\cf}\max\limits_{x\in A}\varphi(x)$ and suppose that
$$\sup\varphi(B)<c.$$
Then, for any sequence of sets $(A_n)_n$ in $\cf$ such that $\lim\limits_{n}\sup\limits_{A_n}\varphi=c$, there exists a sequence $(x_n)_n$ in $X$ such that
\begin{itemize}
	\item [(i)] $\lim\limits_{n}\varphi(x)=c$
	\item [(ii)] $\lim\limits_{n}\|d\varphi(x_n)\|=0$
	\item [(iii)] $\lim\limits_{n} dist(x_n,A_n)=0.$
\end{itemize}
Moreover, if $d\varphi$ is uniformly continuous, then $x_n$ can be chosen to be in $A_n$ for each $n$.
\el
To apply Lemma \ref{le:230820-2} more intuitively, we establish the follwing lemma.
\begin{lemma}\label{20230324_2}
Let $\Upsilon\subset \R \times T(a,b)$ be a closed set and $\tilde I\in C^1(\R\times E,\R)$. For $\theta\in \R$,  we define
$$\tilde{\Gamma}:=\left\{\tilde{\gamma}\in C([0,1],\R\times T(a,b)):\tilde{\gamma}(0)\in\Upsilon,\tilde{I}(\tilde{\gamma}(1))<\theta\right\}.$$
If $\tilde{\Gamma}\neq \emptyset$ and $\tilde{I}$ satisfies
$$c:=\inf\limits_{\tilde{\gamma}\in \tilde{\Gamma}}\max\limits_{t\in[0,1]}\tilde{I}(\tilde{\gamma}(t))>d:=\sup\limits_{\tilde{\gamma}\in \tilde{\Gamma}}\max\{\tilde{I}(\tilde{\gamma}(0)),\tilde{I}(\tilde{\gamma}(1))\},$$
then for any $\varepsilon\in\left(0,\frac{c-d}{2}\right),\delta>0$ and $\tilde{\gamma}_*\in\tilde{\Gamma}$ such that
\beq
\sup\limits_{t\in[0,1]}\tilde{I}(\tilde{\gamma}_*(t))\leq c+\varepsilon,
\eeq
there exists $(\tau,\phi)\in \R\times T(a,b)$ such that
\begin{itemize}
\item[(i)] $c-2\varepsilon\leq\tilde{I}(\tau,\phi)\leq c+2\varepsilon;$
\item[(ii)]$\min\limits_{t\in[0,1]}\|(\tau,\phi)-\tilde{\gamma}_*(t)\|_{\R\times E}\leq 2\delta;$
\item[(iii)]$\|\tilde{I}|'_{\R\times T(a,b)}(\tau,\phi)\|\leq\frac{8\varepsilon}{\delta}.$
\end{itemize}
\end{lemma}
\bp
Since the proof is standard, we will not provide a detailed proof for this. For the scalar case,one can see the details in \cite{Chen2022a}.
\ep
Basing on the above preparations, we are ready to present our main result of this subsection.
\bl\label{le:230821-2}
Under the assumptions of Lemma \ref{le:230820-1}, there exists a nonnegative $(PSP)_{M_{\nu}(a,b)}$ sequence $\{(u_n,v_n)\}$ for $I_{\nu}|'_{T(a,b)}$, i.e, $u_n,v_n\geq0$ with
\beq\label{230820-6}
I_{\nu}(u_n,v_n)\to M_{\nu}(a,b), I_{\nu}|'_{T(a,b)}(u_n,v_n)\to 0 ~\hbox{and}~ P_{\nu}(u_n,v_n)\to0.
\eeq
\el
\bp
Let $\tilde I_{\nu}$ be defined by \eqref{230820-5} and set
$$\tilde\Gamma_{\nu}:=\{\tilde\gamma\in C([0,1],\R\times T(a,b)):\tilde\gamma(0)=(0,u,v),\tilde I_{\nu}(\tilde\gamma(1))<2m_{\nu}(a,b)\}.$$
For any $\gamma=(\gamma_1,\gamma_2)\in\Gamma_{\nu}$, let $\tilde\gamma(t):=(0,\gamma_1(t),\gamma_2(t))$, it is easy to check that $\tilde\gamma\in\tilde\Gamma_{\nu}$. Hence $\tilde\Gamma_{\nu}\neq\emptyset$ and we can define
$$\tilde M_{\nu}(a,b):=\inf\limits_{\tilde\gamma\in\tilde\Gamma_{\nu}}\max\limits_{t\in[0,1]}\tilde I_{\nu}(\tilde\gamma(t)).$$
By the arbitrariness of $\gamma$, we get that $\tilde M_{\nu}(a,b)\leq M_{\nu}(a,b)$. Now we claim that $M_{\nu}(a,b)=\tilde M_{\nu}(a,b)$. Indeed, on the other hand, for any $\tilde\gamma\in\tilde\Gamma_{\nu}$, we denote it by $\tilde\gamma(t)=(\gamma_0(t),\gamma_1(t),\gamma_2(t))$, it is trivial that $\gamma_0(0)=0$. We define $\gamma(t)=(\Theta(\gamma_1(t),\gamma_0(t)),\Theta(\gamma_2(t),\gamma_0(t)))$, then we can check that $\gamma(0)=(\gamma_1(0),\gamma_2(0))=(u,v)$ and
$\tilde I_{\nu}(\tilde\gamma(t))=I_{\nu}(\Theta(\gamma_1(t),\gamma_0(t)),\Theta(\gamma_2(t),\gamma_0(t)))=I_{\nu}(\gamma(t))$, hence $\gamma\in\Gamma_{\nu}$. Then we get that $M_{\nu}(a,b)\leq\tilde M_{\nu}(a,b)$ and thus $M_{\nu}(a,b)=\tilde M_{\nu}(a,b)$. Moreover, combining with Lemma \ref{le:230820-1}, we conclude that
\beq\label{230821-1}
\tilde M_{\nu}(a,b):=
\inf\limits_{\tilde\gamma\in\tilde\Gamma_{\nu}}\max\limits_{t\in[0,1]}\tilde I_{\nu}(\tilde\gamma(t))\geq k_0>0\geq\sup\limits_{\tilde\gamma\in\tilde\Gamma_{\nu}}\max\{\tilde I_{\nu}(\tilde\gamma(0)),\tilde I_{\nu}(\tilde\gamma(1))\}.
\eeq
For any $n\in\N$, there exists $\gamma_n\in\Gamma_{\nu}$, $\gamma_n(t)=(\gamma_{n,1}(t),\gamma_{n,2}(t))$ such that $M_{\nu}(a,b)\leq\max\limits_{t\in[0,1]}I_{\nu}(\gamma_n(t))\leq M_{\nu}(a,b)+\frac{1}{n}$. One can check that $|\gamma_n|\in\Gamma_{\nu}$ since $I_{\nu}$ is even, where $|\gamma_n|(t):=(|\gamma_{n,1}(t)|,|\gamma_{n,2}(t)|)$. So without loss of generality, we may assume that $\gamma_n(t):=|\gamma_n|(t)$. Let $\tilde\gamma_n(t)=(0,\gamma_n(t))$, then $\tilde I_{\nu}(\tilde\gamma_n(t))=I_{\nu}(\gamma_n(t))$ and thus
\beq\label{230821-2}
\max\limits_{t\in[0,1]}\tilde I_{\nu}(\tilde\gamma_n(t))\leq\tilde M_{\nu}(a,b)+\frac{1}{n}.
\eeq
Basing on \eqref{230821-1} and \eqref{230821-2}, by Lemma \ref{20230324_2} we know that there exists a $(PS)_{\tilde M_{\nu}(a,b)}$ sequence $(\tilde s_n,\tilde u_n,\tilde v_n)$ for $\tilde I_{\nu}|_{\R\times T(a,b)}$ satisfying
\begin{itemize}
	\item[(i)] $\tilde M_{\nu}(a,b)-\frac{2}{n}\leq\tilde{I_{\nu}}(\tilde s_n,\tilde u_n,\tilde v_n)\leq \tilde M_{\nu}(a,b)+\frac{2}{n};$
	\item[(ii)]$\min\limits_{t\in[0,1]}\|(\tilde s_n,\tilde u_n,\tilde v_n)-(0,\gamma_n(t))\|_{\R\times E}\leq \frac{2}{\sqrt{n}};$
	\item[(iii)]$\|\tilde{I_{\nu}}|'_{\R\times T(a,b)}(\tilde s_n,\tilde u_n,\tilde v_n)\|\leq\frac{8}{\sqrt{n}}.$
\end{itemize}
Let $(u_n,v_n)=(\Theta(\tilde u_n,\tilde s_n),\Theta(\tilde v_n,\tilde s_n))$, then
$$I_{\nu}(u_n,v_n)=\tilde I_{\nu}(\tilde s_n,\tilde u_n,\tilde v_n)\to \tilde M_{\nu}(a,b)=M_{\nu}(a,b).$$
Combining with \eqref{230820-3} and \eqref{230820-4},  we obtain
$$P_{\nu}(u_n,v_n)=P_{\nu}(\Theta(\tilde u_n,\tilde s_n),\Theta(\tilde v_n,\tilde s_n))\to 0$$
and
$$I_\nu|'_{T(a,b)}(u_n,v_n)=I_\nu|'_{T(a,b)}(\Theta(\tilde u_n,\tilde s_n),\Theta(\tilde v_n,\tilde s_n))\to0.$$
In conclusion, \eqref{230820-6} has been obtained. By (ii) we conclude that $(\tilde u_n^-,\tilde v_n^-)\rightarrow (0,0)$ in $E$. So $\{(u_n^+, v_n^+)\}$ is also a  $(PSP)_{M_\nu(a,b)}$ sequence. Hence, without loss of generality, we can replace $(u_n,v_n)$ by $(u_n^+, v_n^+)$ and conclude the conclusions.
\ep

\subsection{Estimation of the mountain pass level }
Since the nonlinearity has a Sobolev critical growth, one suitable estimation on $M_{\nu}(a,b)$ is needed to recover the compactness of $(PSP)_{M_{\nu}(a,b)}$ sequence $\{(u_n,v_n)\}$ in $E_{rad}$.  This is also the key to solve the open problem raised by Bartsch et al.\cite[Remark 1.3-a)]{Bartsch2023}.

Let $A_N:=[N(N-2)]^{\frac{N-2}{4}}$ and define $U_n(x):=\Theta_n(|x|)\in H_{rad}^{1}(\R^N)$, where
\beq
\Theta_n(r)=A_N\left\{\begin{array}{lcl}
\left(\frac{n}{1+n^{2}r^{2}}\right)^{\frac{N-2}{2}}, &0\leq r<1;\\
\left(\frac{n}{1+n^2}\right)^{\frac{N-2}{2}}(2-r), &1\leq r<2;\\
0, &r\geq 2.
\end{array}
\right.
\eeq
By a direct computation, the following estimations hold, some can be seen in \cite{Chen2022a}:
\beq
\begin{aligned}
\|U_n\|_2^2=O(\frac{\xi(n)}{n^2}), \quad ~\hbox{as}~ n\to\infty,
\end{aligned}
\eeq
where
\beq\label{0427_2}
\xi(n):=\int_0^n\frac{s^{N-1}}{(1+s^2)^{N-2}}ds=
\begin{cases}
O(n),\quad\quad &\hbox{if}~N=3;\\
O\big(\ln(1+n^2)\big), \quad &\hbox{if}~N=4.
\end{cases}
\eeq
\beq\label{0427_3}
\begin{aligned}
\|\nabla U_n\|_2^2&={\int_{\R^N}|\nabla U_n|^2\ud x}=S^{\frac{N}{2}}+O\left(\frac{1}{n^{N-2}}\right), ~\hbox{as}~ n\to \infty.
\end{aligned}
\eeq

\beq\label{0427_4}
\begin{aligned} \|U_n\|^{2^*}_{2^*}=S^{\frac{N}{2}}+O\left(\frac{1}{n^N}\right), ~\hbox{as}~ n\to \infty.
\end{aligned}
\eeq

Let $(u,v)$ be the local minimizer given in Theorem \ref{th:230820-1}. By elliptic regularity result, we remark that $u$ and $v$ are positive radial decreasing smooth functions. So, we also have the following estimations (see \cite{Vicentiu2023})
\beq
\int_{\R^N}uU_n\ud x=O\left(\frac{1}{n^{\frac{N-2}{2}}}\right), ~\hbox{as}~ n\to \infty.
\eeq
\beq
\int_{\R^N}uU_n^{2^*-1}\ud x=O\left(\frac{1}{n^{\frac{N-2}{2}}}\right), ~\hbox{as}~ n\to \infty
\eeq
Notice that
\beq\label{230814-1}
\|\nabla u\|_2^2+\lambda_1\|u\|_2^2=\mu_1\|u\|_{2^*}^{2^*}+\nu\alpha\int_{\R^N}|u|^{\alpha}|v|^{\beta}\ud x,
\eeq
\beq\label{230903-2}
\|\nabla v\|_2^2+\lambda_2\|v\|_2^2=\mu_2\|v\|_{2^*}^{2^*}+\nu\beta\int_{\R^N}|u|^{\alpha}|v|^{\beta}\ud x,
\eeq
and
\beq\label{230814-2}
P_{\nu}(u,v)=\|\nabla u\|_2^2+\|\nabla v\|_2^2-\mu_1\|u\|_{2^*}^{2^*}-\mu_2\|v\|_{2^*}^{2^*}-\nu \gamma_{\alpha+\beta}\int_{\R^N}|u|^{\alpha}|v|^{\beta}\ud x=0.
\eeq
Hence, combining with \eqref{230814-1},\eqref{230903-2} and \eqref{230814-2}, we conclude that
\beq\label{230815-1}
\lambda_1=\frac{1}{a^2}\left(\|\nabla v\|_2^2-\mu_2\|v\|_{2^*}^{2^*}-(\gamma_{\alpha+\beta}-\alpha)\nu\int_{\R^N}u^{\alpha}v^{\beta}\ud x\right)
\eeq
\beq\label{230903-3}
\lambda_2=\frac{1}{b^2}\left(\|\nabla u\|_2^2-\mu_1\|u\|_{2^*}^{2^*}-(\gamma_{\alpha+\beta}-\beta)\nu\int_{\R^N}u^{\alpha}v^{\beta}\ud x\right)
\eeq

Multiplying the first equation in \eqref{eq:230927-e1} by $U_n$, we have that
\beq\label{230815-7}
\int_{\R^N}\nabla u\nabla U_n\ud x+\lambda_1\int_{\R^N}uU_n\ud x=\mu_1\int_{\R^N}u^{2^*-1}U_n\ud x+\nu\alpha\int_{\R^N}u^{\alpha-1}v^{\beta}U_n\ud x.
\eeq
Let $\tau=\frac{\|u+tU_n\|_2}{a}$, then
\beq\label{0428_4}
\tau^2=1+\frac{2t}{a^2}\int_{\R^N}u U_n\ud x+t^2\left[O(\frac{\xi(n)}{n^2})\right],\quad ~\hbox{as}~ n\to\infty,
\eeq
see \cite[formula(4.41)]{Vicentiu2023}.
We remark that $\tau>1$ depends on $n$ and $t$, it can also be denoted by $\tau_{n,t}$. In particular, $\tau\rightarrow 1$ as $t\rightarrow 0$ uniformly for $n$ and $\tau\rightarrow +\infty$ as $t\rightarrow +\infty$ for any $n\in \N$ fixed. Furthermore, $\tau\rightarrow 1$ as $n\rightarrow +\infty$ uniformly for bounded $t$.
Define
$$W_{n,t}=\tau\star\frac{u+tU_n}{\tau},$$
then
$$\|W_{n,t}\|_2^2=\|\frac{u+tU_n}{\tau}\|_2^2=a^2, \|\nabla W_{n,t}\|_2^2=\|\nabla (u+tU_n)\|_2^2 ~\hbox{and}~ \|W_{n,t}\|_{2^*}^{2^*}=\|u+tU_n\|_{2^*}^{2^*}.$$
Thus, combining with \eqref{0428_4}, we see that $\{(W_{n,t},\tau\star v):t\geq0\}$ is a curve in $T(a,b)$ for any fixed $n\in\N$. Let $H_n:\R\rightarrow \R$ be a map defined by
\beq\label{230821-3}
\begin{aligned}
H_n(t):&=I_{\nu}(W_{n,t},\tau\star v)\\
&=\frac{1}{2}\|\nabla(u+tU_n)\|_2^2+\frac{1}{2}\tau^2\|\nabla v\|_2^2-\frac{\mu_1}{2^*}\|u+tU_n\|_{2^*}^{2^*}-\frac{\mu_2}{2^*}\tau^{2^*}\|v\|_{2^*}^{2^*}\\
&\,\,\,\,\,\,-\nu\tau^{\gamma_{\alpha+\beta}-\alpha}\int_{\R^N}|u+tU_n|^{\alpha}|v|^{\beta}\ud x.
\end{aligned}
\eeq
\bl\label{le:230815-1}
Let $H_n(t)$ be given above, then there exist $N_0\in\N$ and $T>0$ such that $H_n(T)<2m_{\nu}(a,b)<0$ for all $n\geq N_0$. Define $\gamma_n(t):=(W_{n,tT}, \tau_{n,tT}\star v)$, then $\gamma_n\in\Gamma_\nu$. Moreover, there exists $0<t_n<T$ such that $H_n(t_n)=\max\limits_{t>0}H_n(t)$ and $0<\inf\limits_{n\geq N_0}t_n\leq\sup\limits_{n\geq N_0}t_n<+\infty$.
\el
\bp
For any $n\in\N$, $H_n(0)=I_{\nu}(W_{n,0},v)=I_{\nu}(u,v)=m_{\nu}(a,b)<0$. We {\bf claim} that there exists $N_0\in\N$ such that $H_n(t)\to-\infty$ as $t\to+\infty$ uniformly for all $n\geq N_0$. Indeed, by \eqref{230821-3}
we have that
\begin{align*}
H_n(t)&\leq\frac{1}{2}\|\nabla u\|_2^2+\frac{t^2}{2}\|\nabla U_n\|_2^2+t\int_{\R^N} \nabla u \nabla U_n\ud x+\frac{\tau^2}{2}\|\nabla v\|_2^2\\
&\,\,\,\,\,\,\,\,\,\,\,\,\,\,\,\,\,\,-\frac{\mu_1t^{2^*}}{2^*}\|U_n\|_{2^*}^{2^*}-\frac{\mu_2\tau^{2^*}}{2^*}\|v\|_{2^*}^{2^*}-\nu\tau^{\gamma_{\alpha+\beta}-\alpha}\int_{\R^N}t^{\alpha}|U_n|^{\alpha}|v|^{\beta}\ud x\nonumber\\
&\leq\frac{1}{2}\|\nabla u\|_2^2+\frac{t^2}{2}(S^{\frac{N}{2}}+o_n(1))-\frac{\mu_1t^{2^*}}{2^*}(S^{\frac{N}{2}}+o_n(1))+to_n(1)\\
&\,\,\,\,\,\,+\frac{1}{2}\|\nabla v\|_2^2(1+t^2o_n(1)+to_n(1))\\
&=-\frac{\mu_1}{2^*}(S^{\frac{N}{2}}+o_n(1))t^{2^*}+\frac{1}{2}(S^{\frac{N}{2}}+o_n(1))t^2+o_n(1)+\frac{1}{2}(\|\nabla u\|_2^2+\|\nabla v\|_2^2)\\
&\to-\infty ~\hbox{as}~ t\to+\infty ~\hbox{uniformly for large}~ n.
\end{align*}
Hence, we can find $N_0\in\N$ and $T>0$ large enough such that $H_n(T)<2m_{\nu}(a,b)$ for all $n\geq N_0$. Thus $\gamma_n\in\Gamma_\nu$. Combining with Lemma \ref{le:230820-1}, for any $n\geq N_0$, there there exists $0<t_n<T$ such that
$$H_n(t_n)=\max\limits_{t>0}H_n(t)=\max\limits_{t\in[0,1]}I_{\nu}(\gamma_n(t))\geq k_0>0.$$
So, $\sup\limits_{n\geq N_0}t_n\leq T<+\infty$.

On the other hand, we {\bf claim} that $\inf\limits_{n\geq N_0}t_n>0$. If not, there exists a subsequence $n_k\rightarrow +\infty$ such that $t_{n_k}\rightarrow 0$.
Then it is trivial that $W_{n_k,t_{n_k}}\to u$ in $H^1(\R^N)$,$\tau_{n_k,t_{n_k}}\star v\to v$ in $H^1(\R^N)$. Thus, $H_{n_k}(t_{n_k})=I_{\nu}(W_{n_k,t_{n_k}},\tau_{n_k,t_{n_k}}\star v)\to I_{\nu}(u,v)=m_{\nu}(a,b)<0$, which contradicts with $H_n(t_n)\geq k_0>0$ for all $n\geq N_0$.
\ep

By Lemma \ref{le:230815-1} we know that for all $n\geq N_0$, $M_{\nu}(a,b)\leq\max\limits_{t\in[0,1]}I(\gamma_n(t))=H_n(t_n)$ and $t_n$ is bounded. For the sake of convinence, we denote $\tau_{n,t_n}$ by $\tau_n$. By \eqref{0428_4}, we have
$$\tau_n^2=1+\frac{2t_n}{a^2}\int_{\R^N}u U_n\ud x+t_n^2O(\frac{\xi(n)}{n^2}),\quad ~\hbox{as}~ n\to\infty.$$
Then for any $m\in\R$, we also have
\beq\label{230815-2}
\begin{aligned}
\tau_n^m=&\left(1+\frac{2t_n}{a^2}\int_{\R^N}u U_n\ud x+t_n^2O(\frac{\xi(n)}{n^2})\right)^{\frac{m}{2}}\\
=&1+\frac{mt_n}{a^2}\int_{\R^N}u U_n\ud x+O(\frac{\xi(n)}{n^2}).
\end{aligned}
\eeq
Hence for any $m\neq0$, we conclude that
\beq
\tau_n^m-1=O\left(\frac{1}{n^{\frac{N-2}{2}}}\right)~\hbox{as}~ n\to +\infty.
\eeq
\bl\label{le:230815-2}
Suppose \eqref{230815-8}. Let $\nu_0$ be given by Theorem \ref{th:230820-1}-(ii). $M_{\nu}(a,b), m_{\nu}(a,b)$ are defined by \eqref{230902-2} and \eqref{230921-1} respectively. Then for any $\nu\in (0,\nu_0)$, it follows that
\beq\label{230903-6}
M_{\nu}(a,b)<m_{\nu}(a,b)+\frac{1}{N}\min\{{\mu_1}^{\frac{2-N}{2}},{\mu_2}^{\frac{2-N}{2}}\}S^{\frac{N}{2}}.
\eeq
\el
\bp
For any $t>0$,  we remark that
\beq\label{230815-6}
(1+t)^p\geq
\begin{cases}
1+t^p+pt^{p-1}+pt, &p\geq3;\\
1+pt, &p\geq1.	
\end{cases}
\eeq

Without loss of generality, we may assume that $\mu_1\geq \mu_2$, the case of $\mu_1<\mu_2$ can be proved by a similar way.

Noting that $2^*\geq3$ for $N\in\{3,4\}$, we have
{\allowdisplaybreaks
\begin{align}\label{230815-3}
H_n(t_n)=&\frac{1}{2}\|\nabla(u+t_nU_n)\|_2^2+\frac{1}{2}\tau_n^2\|\nabla v\|_2^2-\frac{\mu_1}{2^*}\|u+t_nU_n\|_{2^*}^{2^*}-\frac{\mu_2}{2^*}\tau_n^{2^*}\|v\|_{2^*}^{2^*}\nonumber\\
&-\nu\tau_{n}^{\gamma_{\alpha+\beta}-\alpha}\int_{\R^N}|u+t_nU_n|^{\alpha}|v|^{\beta}\ud x\nonumber\\
\hbox{by \eqref{230815-6} }\leq&\frac{1}{2}\left(\|\nabla u\|_2^2+t_n^2\|\nabla U_n\|_2^2+2t_n\int_{\R^N}\nabla u\nabla U_n\ud x\right)+\frac{1}{2}\tau_n^2\|\nabla v\|_2^2\nonumber\\
&-\frac{\mu_1}{2^*}\left(\|u\|_{2^*}^{2^*}+t_{n}^{2^*}\|U_n\|_{2^*}^{2^*}+2^* t_n\int_{\R^N}u^{2^*-1}U_n\ud x+2^*t_{n}^{2^*-1}\int_{\R^N}u U_{n}^{2^*-1}\ud x\right)\nonumber\\
&-\frac{\mu_2}{2^*}\tau_{n}^{2^*}\|v\|_{2^*}^{2^*}-\nu\tau_{n}^{\gamma_{\alpha+\beta}-\alpha}\int_{\R^N}|u|^{\alpha}|v|^{\beta}\ud x-\nu\alpha\tau_{n}^{\gamma_{\alpha+\beta}-\alpha}  t_n\int_{\R^N}u^{\alpha-1}v^{\beta}U_n\nonumber\\
\hbox{by \eqref{230815-7}} =&\frac{1}{2}\|\nabla u\|_2^2+\frac{1}{2}\tau_n^2\|\nabla v\|_2^2-\frac{\mu_1}{2^*}\|u\|_{2^*}^{2^*}-\frac{\mu_2}{2^*}\tau_{n}^{2^*}\|v\|_{2^*}^{2^*}-\nu\int_{\R^N}|u|^{\alpha}|v|^{\beta}\ud x+\frac{t_n^2}{2}\|\nabla U_n\|_2^2\nonumber\\
&-\frac{\mu_1t_{n}^{2^*}}{2^*}\|U_n\|_{2^*}^{2^*}-\lambda_1t_n\int_{\R^N}uU_n\ud x+\nu\alpha(1-\tau_{n}^{\gamma_{\alpha+\beta}-\alpha}) t_n\int_{\R^N}u^{\alpha-1}v^{\beta}U_n\ud x\nonumber\\
&-\mu_1t_{n}^{2^*-1}\int_{\R^N}u U_{n}^{2^*-1}\ud x+\nu(1-\tau_{n}^{\gamma_{\alpha+\beta}-\alpha})\int_{\R^N}|u|^{\alpha}|v|^{\beta}\ud x\nonumber\\
\end{align}
}
Recalling \eqref{230815-1} and \eqref{230815-2},
{\allowdisplaybreaks
\begin{align}\label{230815-4}
&\frac{1}{2}\tau_n^2\|\nabla v\|_2^2-\frac{\mu_2}{2^*}\tau_{n}^{2^*}\|v\|_{2^*}^{2^*}-\lambda_1t_n\int_{\R^N}uU_n\ud x+\nu(1-\tau_{n}^{\gamma_{\alpha+\beta}-\alpha})\int_{\R^N}|u|^{\alpha}|v|^{\beta}\ud x\nonumber\\	
=&\frac{1}{2}\left(1+\frac{2t_n}{a^2}\int_{\R^N}u U_n\ud x+O(\frac{\xi(n)}{n^2})\right)\|\nabla v\|_2^2-\frac{\mu_2}{2^*}\left(1+\frac{2^*t_n}{a^2}\int_{\R^N}u U_n\ud x+O(\frac{\xi(n)}{n^2})\right)\|v\|_{2^*}^{2^*}\nonumber\\
&-\frac{t_n}{a^2}\left(\|\nabla v\|_2^2-\mu_2\|v\|_{2^*}^{2^*}-(\gamma_{\alpha+\beta}-\alpha)\nu\int_{\R^N}u^{\alpha}v^{\beta}\ud x\right)\int_{\R^N}uU_n\ud x\nonumber\\
&-\nu\left(\frac{(\gamma_{\alpha+\beta}-\alpha)t_n}{a^2}\int_{\R^N}uU_n\ud x+O(\frac{\xi(n)}{n^2})\right)\int_{\R^N}u^{\alpha}v^{\beta}\ud x\nonumber\\
=&\frac{1}{2}\|\nabla v\|_2^2-\frac{\mu_2}{2^*}\|v\|_{2^*}^{2^*}+O(\frac{\xi(n)}{n^2})
\end{align}
}
Hence, combining \eqref{230815-3} and \eqref{230815-4} we have that
{\allowdisplaybreaks
\begin{align}\label{230815-5}
H_n(t_n)&\leq\frac{1}{2}\|\nabla u\|_2^2+\frac{1}{2}\|\nabla v\|_2^2-\frac{\mu_1}{2^*}\|u\|_{2^*}^{2^*}-\frac{\mu_2}{2^*}\|v\|_{2^*}^{2^*}-\nu\int_{\R^N}|u|^{\alpha}|v|^{\beta}\ud x\nonumber\\
&\,\,\,\,\,+\frac{t_n^2}{2}\|\nabla U_n\|_2^2-\frac{\mu_1t_{n}^{2^*}}{2^*}\|U_n\|_{2^*}^{2^*}-\mu_1t_{n}^{2^*-1}\int_{\R^N}u U_{n}^{2^*-1}\ud x\nonumber\\
&\,\,\,\,\,+\nu\alpha(1-\tau_{n}^{\gamma_{\alpha+\beta}-\alpha}) t_n\int_{\R^N}u^{\alpha-1}v^{\beta}U_n\ud x+O(\frac{\xi(n)}{n^2})\nonumber\\
=&m_{\nu}(a,b)+\frac{t_n^2}{2}\left(S^{\frac{N}{2}}+O\left(\frac{1}{n^{N-2}}\right)\right)-\frac{t_{n}^{2^*}}{2^*}\left(\mu_1S^{\frac{N}{2}}+O\left(\frac{1}{n^N}\right)\right)\nonumber\\
&-t_{n}^{2^*-1}O\left(\frac{1}{n^{\frac{N-2}{2}}}\right)-t_nO\left(\frac{1}{n^{\frac{N-2}{2}}}\right)O\left(\frac{1}{n^{\frac{N-2}{2}}}\right)+O(\frac{\xi(n)}{n^2})\nonumber\\
=&m_{\nu}(a,b)+(\frac{t_n^2}{2}-\frac{\mu_1t_{n}^{2^*}}{2^*})S^{\frac{N}{2}}-O\left(\frac{1}{n^{\frac{N-2}{2}}}\right)+O\left(\frac{1}{n^{N-2}}\right)\\
&-O\left(\frac{1}{n^N}\right)+O\left(\frac{\xi(n)}{n^2}\right)\nonumber
\end{align}
}
Noting that $\frac{t^2}{2}-\frac{\mu_1t^{2^*}}{2^*}\leq\frac{1}{N}{\mu_1}^{\frac{2-N}{2}}$ for all $t>0$. By \eqref{230815-5}, there exists $N_1\in\N$ such that
$$
-O\left(\frac{1}{n^{\frac{N-2}{2}}}\right)+O\left(\frac{1}{n^{N-2}}\right)-O\left(\frac{1}{n^N}\right)+O\left(\frac{\xi(n)}{n^2}\right)<0,\quad \forall n\geq N_1,
$$
then $H_n(t_n)<m_{\nu}(a,b)+\frac{1}{N}{\mu_1}^{\frac{2-N}{2}}S^{\frac{N}{2}}$ for any $n\geq\max\{N_0,N_1\}$ and thus
\beq\label{230903-4}
M_{\nu}(a,b)<m_{\nu}(a,b)+\frac{1}{N}{\mu_1}^{\frac{2-N}{2}}S^{\frac{N}{2}}.
\eeq
The proof is finished.
\ep

\subsection{Compactness and the proof of Theorem \ref{th:230821-1}}

Now, we are ready to prove our main result.

\noindent
{\bf Proof of Theorem \ref{th:230821-1}:}
By Lemma \ref{le:230821-2}, there exists a nonnegative $(PSP)_{M_{\nu}(a,b)}$ sequence $\{(u_n,v_n)\}\subset T(a,b)$ for $I_{\nu}|'_{T(a,b)}$, i.e. there exist $\{\lambda_{1,n}\}, \{\lambda_{2,n}\}\subset\R$ such that $\{(\lambda_{1,n}, \lambda_{2,n},u_n,v_n)\}$ satisfies
\beq\label{230904-3}
\begin{cases}
-\Delta u_n+\lambda_{1,n}u_n=\mu_1 u_{n}^{2^*-1}+\nu\alpha u_{n}^{\alpha-1}v_{n}^{\beta}+o_n(1)~\hbox{in}~H^{-1}(\R^N),\\
-\Delta v_n+\lambda_{2,n}v_n=\mu_2 v_{n}^{2^*-1}+\nu\beta u_{n}^{\alpha}v_{n}^{\beta-1}+o_n(1)~\hbox{in}~H^{-1}(\R^N),
\end{cases}
\eeq
and $I_{\nu}(u_n,v_n)\to M_{\nu}(a,b)$, $P_{\nu}(u_n,v_n)\to0$.
Hence, recalling \eqref{230904-1} and \eqref{230904-2}, by \eqref{230901-1} we have
\begin{align*}
&M_{\nu}(a,b)+o_n(1)=I_{\nu}(u_n,v_n)-\frac{1}{2^*}P_{\nu}(u_n,v_n)\\
=&(\frac{1}{2}-\frac{1}{2^*})(\|\nabla u_n\|_2^2+\|\nabla v_n\|_2^2)-\nu\frac{2^*-\gamma_{\alpha+\beta}}{2^*}\int_{\R^N}|u_n|^{\alpha}|v_n|^{\beta}\ud x\\
\geq&(\frac{1}{2}-\frac{1}{2^*})(\|\nabla u_n\|_2^2+\|\nabla v_n\|_2^2)-\nu\frac{2^*-\gamma_{\alpha+\beta}}{2^*}C(N,\alpha,\beta)(a^2+b^2)^{\frac{\alpha+\beta-\gamma_{\alpha+\beta}}{2}}(\|\nabla u_n\|_2^2+\|\nabla v_n\|_2^2)^{\frac{\gamma_{\alpha+\beta}}{2}},
\end{align*}
which implies $\|\nabla u_n\|_2^2+\|\nabla v_n\|_2^2<+\infty$ since $0<\frac{\gamma_{\alpha+\beta}}{2}<1$, hence $\{(u_n,v_n)\}$ is bounded in $E_{rad}$. Thus, by \eqref{230904-3}, we obtain
$$\lambda_{1,n}+o_n(1)=\frac{1}{a^2}\left(-\|\nabla u_n\|_2^2+\mu_1\|u_n\|_{2^*}^{2^*}+\nu\alpha\int_{\R^N}u_{n}^{\alpha}v_{n}^{\beta}\ud x\right)$$
and
$$\lambda_{2,n}+o_n(1)=\frac{1}{b^2}\left(-\|\nabla v_n\|_2^2+\mu_2\|v_n\|_{2^*}^{2^*}+\nu\beta\int_{\R^N}u_{n}^{\alpha}v_{n}^{\beta}\ud x\right),$$
which imply that $\{\lambda_{1,n}\}$ and $\{\lambda_{2,n}\}$ are also bounded. Going up to a subsequence, there exists $(\bar \lambda_1,\bar\lambda_2,u_0,v_0)\in\R^2\times E_{rad}$ such that $(u_n,v_n)\rightharpoonup(u_0,v_0)$ in $E_{rad}$ and  $\lambda_{1,n}\to\bar\lambda_1, \lambda_{2,n}\to\bar\lambda_2$. Then one can see that $(u_0,v_0)$ is a nonnegative weak radial solution to
\beq\label{230904-4}
\begin{cases}
-\Delta u+\bar \lambda_1 u=\mu_1 u^{2^*-1}+\nu\alpha u^{\alpha-1}v^{\beta} ~\hbox{in}~ \R^N,\\
-\Delta v+\bar \lambda_2 v=\mu_2 v^{2^*-1}+\nu\beta u^{\alpha}v^{\beta-1} ~\hbox{in}~ \R^N.\\
\end{cases}
\eeq
It follows $P_{\nu}(u_0,v_0)=0$.
Now we claim that $u_0\not\equiv0$ and $v_0\not\equiv0$.

Suppose $u_0\equiv0$, then by $P_{\nu}(u_0,v_0)=0$ we have $\bar\lambda_2=0$. So $v_0$ is a radial solution to
$$-\Delta v=\mu_2 v^{2^*-1}, v\geq0.$$
If $v_0\not\equiv0$, then we know that $v_0=\mu_{2}^{\frac{N-2}{4}}U$, where $U$ is the unique radial solution of \eqref{eq:230921-1}. It is well known that $U\not\in L^2(\R^N)$ for $N\in\{3,4\}$, which contradicts with $v_0\in H^1(\R^N)$. Hence, $v_0\equiv0$. Since $u_n\rightharpoonup u_0\equiv0$ in $H_{rad}^{1}(\R^N)$, then $\int_{\R^N}|u_n|^{\alpha}|v_n|^{\beta}\ud x\to0$ and thus
\beq\label{230904-5}
P_{\nu}(u_n,v_n)=\|\nabla u_n\|_2^2+\|\nabla v_n\|_2^2-\mu_1\|u_n\|_{2^*}^{2^*}-\mu_2\|v_n\|_{2^*}^{2^*}=o_n(1).
\eeq
Assume that $\|\nabla u_n\|_2^2\to h_1^2$ and $\|\nabla v_n\|_2^2\to h_2^2$, then by \eqref{230904-5} we have
\begin{align}\label{230904-6}
M_{\nu}(a,b)+o_n(1)=&I_{\nu}(u_n,v_n)=\frac{1}{2}(\|\nabla u_n\|_2^2+\|\nabla v_n\|_2^2)-\frac{\mu_1}{2^*}\|u_n\|^{2^*}_{2^*}-\frac{\mu_2}{2^*}\|v_n\|^{2^*}_{2^*}+o_n(1)\nonumber\\
=&\frac{1}{N}(\|\nabla u_n\|_2^2+\|\nabla v_n\|_2^2)+o_n(1)
=\frac{1}{N}(h_1^2+h_2^2)+o_n(1).
\end{align}
If $h_1^2+h_2^2=0$, then $I_{\nu}(u_n,v_n)\to0$, which contradicts with $M_{\nu}(a,b)>0$. On the other hand, according to critical Sobolev inequality, it follows from \eqref{230904-5} that
$h_1^2+h_2^2\leq S^{-\frac{2^*}{2}}(\mu_1 h_{1}^{2^*}+\mu_2 h_{2}^{2^*})$.
Hence, by \eqref{230904-6}, we have
\begin{align*}
M_{\nu}(a,b)\geq&\frac{1}{N}\min\{h_1^2+h_2^2:0<h_1^2+h_2^2\leq S^{-\frac{2^*}{2}}(\mu_1 h_{1}^{2^*}+\mu_2 h_{2}^{2^*})\}\\
=&\frac{1}{N}\min\{r^2:0<r^2\leq S^{-\frac{2^*}{2}}r^{2^*}(\mu_1 cos^{2^*}\theta+\mu_2 sin^{2^*}\theta), \theta\in [0,\frac{\pi}{2}]\}\\
\geq&\frac{1}{N}\min\{{\mu_1}^{\frac{2-N}{2}},{\mu_2}^{\frac{2-N}{2}}\}S^{\frac{N}{2}}>m_{\nu}(a,b)+\frac{1}{N}\min\{{\mu_1}^{\frac{2-N}{2}},{\mu_2}^{\frac{2-N}{2}}\}S^{\frac{N}{2}},
\end{align*}
which contradicts with \eqref{230903-6}. Thus $u_0\not\equiv0$ and $v_0\not\equiv0$.

Since $\{(u_n,v_n)\}$ is bounded in $E_{rad}$, by Br\'ezis-Lieb Lemma, we have
$$I_{\nu}(u_n,v_n)=I_{\nu}(u_0,v_0)+\frac{1}{2}(\|\nabla (u_n-u_0)\|_2^2+\|\nabla (v_n-v_0)\|_2^2)-\frac{\mu_1}{2^*}\|u_n-u_0\|_{2^*}^{2^*}-\frac{\mu_2}{2^*}\|v_n-v_0\|_{2^*}^{2^*}+o_n(1)$$
and
$$P_{\nu}(u_n,v_n)=P_{\nu}(u_0,v_0)+\|\nabla (u_n-u_0)\|_2^2+\|\nabla (v_n-v_0)\|_2^2-\mu_1 \|u_n-u_0\|_{2^*}^{2^*}-\mu_2\|v_n-v_0\|_{2^*}^{2^*}+o_n(1).$$
Combining with critical Sobolev inequality and by a similar argument above, we have
$\|\nabla (u_n-u_0)\|_2^2+\|\nabla (v_n-v_0)\|_2^2\to0$ or $\|\nabla (u_n-u_0)\|_2^2+\|\nabla (v_n-v_0)\|_2^2\geq \min\{{\mu_1}^{\frac{2-N}{2}},{\mu_2}^{\frac{2-N}{2}}\}S^{\frac{N}{2}}+o_n(1)$. If the latter holds, by Lemma \ref{le:230821-1}-(ii), then we have
\begin{align*}
M_{\nu}(a,b)+o_n(1)=&I_{\nu}(u_n,v_n)\\
=& I_{\nu}(u_0,v_0)+(\frac{1}{2}-\frac{1}{2^*})(\|\nabla (u_n-u_0)\|_2^2+\|\nabla (v_n-v_0)\|_2^2)+o_n(1)\\
\geq&m_{\nu}(\|u_0\|_2,\|v_0\|_2)+\frac{1}{N}\min\{{\mu_1}^{\frac{2-N}{2}},{\mu_2}^{\frac{2-N}{2}}\}S^{\frac{N}{2}}+o_n(1)\\
\hbox{(by Lemma \ref{le:230821-1}-(ii))}\geq& m_{\nu}(a,b)+\frac{1}{N}\min\{{\mu_1}^{\frac{2-N}{2}},{\mu_2}^{\frac{2-N}{2}}\}S^{\frac{N}{2}}+o_n(1),
\end{align*}
which contradicts with \eqref{230903-6}. Hence, $\|\nabla (u_n-u_0)\|_2^2+\|\nabla (v_n-v_0)\|_2^2\to0$, i.e. $(u_n,v_n)\to (u_0,v_0)$ in $D_{0}^{1,2}(\R^N)\times D_{0}^{1,2}(\R^N)$.

Now we {\bf claim} that $\bar\lambda_1>0$ and $\bar\lambda_2>0$. If not, suppose $\bar\lambda_1\leq0$, then
$$
\begin{cases}
-\Delta u_0=-\bar\lambda_1 u_0+\mu_1 u_{0}^{2^*-1}+\nu\alpha u_{0}^{\alpha-1}v_{0}^{\beta}\geq0 ~\hbox{in}~ \R^N,\\
-\Delta v_0=-\bar\lambda_2 v_0+\mu_2 v_{0}^{2^*-1}+\nu\beta u_{0}^{\alpha}v_{0}^{\beta-1} ~\hbox{in}~ \R^N, \\
u_0\geq0, v_0\geq0.
\end{cases}
$$
By the standard ellptic regularity, $u_0$ and $v_0$ are smooth. Hence, by Liouville theorem \cite[Lemma A.2]{Ikoma2014} we get $u_0\equiv0$, this is a contradiction. Similarly, we can show that $\bar\lambda_2>0$.

According to Br\'ezis-Lieb Lemma,
we have $$(I_{\nu}'(u_n,v_n),(u_n,v_n))=(I_{\nu}'(u_0,v_0),(u_0,v_0))+\bar\lambda_1\|u_n-u_0\|_2^2+\bar\lambda_2\|v_n-v_0\|_2^2+o_n(1).$$
So $(u_n,v_n)\to (u_0,v_0)$ in $H$.
And thus $(u_n,v_n)\to(u_0,v_0)\in E_{rad}$.
\hfill$\Box$

\s{The asymptotic behavior and the proof of theorem \ref{th:230913-1}}\label{sec:4}
For simplicity, we write $A\sim B$ if there exist $C_1,C_2>0$ such that $C_1B\leq A\leq C_2B$. Under the condition \eqref{230815-8}, let $(u_{\nu},v_{\nu})$ be the normalized ground state given by Theorem \ref{th:230820-1}, then it is positive, radial symmetric and decreasing. Although here we consider a much more general model than that in \cite{Bartsch2023}, one can obtain the following Lemma \ref{le:230921-1} and Lemma \ref{le:230914-1}, just by a slight modification of \cite[Lemma 5.1]{Bartsch2023} and \cite[Lemma 5.2]{Bartsch2023} respectively. So we prefer to state out the conclusion as follows without a proof.

\bl\label{le:230921-1}
Suppose \eqref{230815-8}, $\nu<\nu_0$ and $\{(u_{\nu},v_{\nu})\}$ are a family of normalized ground state solutions given by Theorem \ref{th:230820-1}, $\{\lambda_{1,\nu},\lambda_{2,\nu}\}$ are corresponding Lagrange multipliers. Then
$$\lambda_{1,\nu}+\lambda_{2,\nu}\sim\|\nabla u_{\nu}\|_2^2+\|\nabla v_{\nu}\|_2^2\sim\nu^{\frac{2}{2-\gamma_{\alpha+\beta}}}, ~\hbox{as}~ \nu\to0^+.$$
\el

\bl\label{le:230914-1}
Under the assumptions of Theorem \ref{th:230820-1}, then there exists $t_{\nu}\in\R$ such that
$$dist_{E}\left(t_{\nu}\star (u_{\nu},v_{\nu}),G(a,b)\right)\to0, ~\hbox{as}~ \nu\to 0^+.$$
Moreover, $t_{\nu}\sim\nu^{\frac{1}{2-\gamma_{\alpha+\beta}}}$ as $\nu\to0^+$.
\el

The asymptotic behavior of the normalized ground state can be studied by a way of \cite[Theorem 1.7a)]{Bartsch2023}. By Lemma \ref{le:230921-1} we see that $\lambda_{1,\nu}+\lambda_{2,\nu}\to0$ as $\nu\to0$, then we conclude that $\lambda_{1,\nu}\to0$ and $\lambda_{2,\nu}\to0$ due to the fact both $\lambda_{1,\nu}$ and $\lambda_{2,v}$ are positive. Basing on this, we can obtain something more as follows.
\bl\label{le:230914-2}
Under the assumptions of Theorem \ref{th:230820-1}-(ii), $\{(u_{\nu},v_{\nu})\}$ is a family of normalized ground state solutions given by Theorem \ref{th:230820-1}-(ii). Then $m_{\nu}(a,b)\to0$ and $\|u_{\nu}\|_{\infty}+\|v_{\nu}\|_{\infty}\to0$ as $\nu\to0^+$.
\el
\bp
By Remark \ref{br:230913-1} or Lemma \ref{le:230921-1}, we know that $(u_{\nu},v_{\nu})\to (0,0)$ in $D_{0}^{1,2}(\R^N)\times D_{0}^{1,2}(\R^N)$. Combining with $(u_{\nu},v_{\nu})\in T(a,b)$, it is trivial that $(u_{\nu},v_{\nu})$ is bounded in $E$. Hence, $m_{\nu}(a,b)=\frac{1}{2}(\|\nabla u_{\nu}\|_2^2+\|\nabla v_{\nu}\|_2^2)-\frac{\mu_1}{2^*}\|u_{\nu}\|^{2^*}_{2^*}-\frac{\mu_2}{2^*}\|v_{\nu}\|^{2^*}_{2^*}-\nu\int_{\R^N}u_{\nu}^{\alpha}v_{\nu}^{\beta}\ud x\to0$ as $\nu\to0$.
Firstly, we {\bf claim} that:
$$\hbox{ there exists $M>0$ such that $\lim\limits_{\nu\to0}(\|u_{\nu}\|_{\infty}+\|v_{\nu}\|_{\infty})\leq M$.}$$
Since $(u_{\nu},v_{\nu})$ is radial symmetric and decreasing, then $\|u_{\nu}\|_{\infty}=u_{\nu}(0)$, $\|v_{\nu}\|_{\infty}=v_{\nu}(0)$. If the claim is false, up to a subsequence, without loss of generality, we can assume that $v_{\nu}(0)\leq u_{\nu}(0)=:L_{\nu}\to+\infty$ as $\nu\to0$. Let
$$\tilde u_{\nu}(y):=\frac{1}{L_{\nu}}u_{\nu}(\frac{y}{L_{\nu}^{\frac{2^*-2}{2}}}),\quad \tilde v_{\nu}(y):=\frac{1}{L_{\nu}}v_{\nu}(\frac{y}{L_{\nu}^{\frac{2^*-2}{2}}}),$$
then $\|\tilde u_{\nu}\|_{\infty}=1$ and $\|\tilde v_{\nu}\|_{\infty}\leq1$. Since $(u_{\nu},v_{\nu})$ is a solution to system  \eqref{eq:230927-e1} with $(\lambda_1,\lambda_2)=(\lambda_{1,\nu}, \lambda_{2,\nu})$, one can check that $\tilde u_{\nu}$ satisfies
\beq
-\Delta\tilde u_{\nu}=-\lambda_{1,\nu}L_{\nu}^{2-2^*}\tilde u_{\nu}+\mu_1 \tilde u_{\nu}^{2^*-1}+\nu\alpha L_{\nu}^{\alpha+\beta-2^*}\tilde{u}_{\nu}^{\alpha-1}\tilde{v}_{\nu}^{\beta}.
\eeq
By $\tilde u_{\nu}\in L^{\infty}(\R^N)$, $L_{\nu}\to+\infty$ and $\lambda_{1,\nu}\to0$, we see that $-\lambda_{1,\nu}L_{\nu}^{2-2^*}\tilde u_{\nu}+\mu_1 \tilde u_{\nu}^{2^*-1}+\nu\alpha L_{\nu}^{\alpha+\beta-2^*}\tilde{u}_{\nu}^{\alpha-1}\tilde{v}_{\nu}^{\beta}\in L^{\infty}(\R^N)$. They by a standard elliptic estimation, going up to a subsequence if necessary, we may assume that $\tilde u_{\nu}\to u_0$ in $C_{loc}^{2}(\R^N)$ as $\nu\rightarrow 0^+$. So $u_0(0)=\lim\limits_{\nu\to0}\tilde u_{\nu}(0)=1$ and thus $u_0\not\equiv0$. In particular, $u_0\geq 0$ is a radial solution to the following equation
\beq\label{230913-2}
-\Delta u=\mu_1 u^{2^*-1}.
\eeq
By the well known classify result, we obtain that $u_0=\mu_1^{\frac{2-N}{4}}U$, where $U$ is the unique positive radial solution of \eqref{eq:230921-1}. Noting that $\lim\limits_{\nu\to0}\|\tilde u_{\nu}\|_{2^*}^{2^*}\geq\|u_0\|_{2^*}^{2^*}=\|\mu_1^{\frac{2-N}{4}}U\|_{2^*}^{2^*}=\mu_1^{-\frac{N}{2}}S^{\frac{N}{2}}$, but by Remark \ref{br:230913-1} we know that $\|\tilde u_{\nu}\|_{2^*}^{2^*}=\|u_{\nu}\|_{2^*}^{2^*}\to0$, this is a contradiction. Hence the claim is true.
Combining with $\lambda_{1,\nu}\to0$, by a standard elliptic regularity again, we know that $u_{\nu}\to u$ in $C_{loc}^{2}(\R^N)$ and $u$ satisfies \eqref{230913-2}. If $u\not\equiv0$, then we can deduce a contradiction similar to the argument above. Hence $u\equiv0$ and thus $\|u_{\nu}\|_{\infty}=u_{\nu}(0)\to0$, i.e., $\|u_{\nu}\|_{\infty}+\|v_{\nu}\|_{\infty}\to0$.
\ep

Our main effort in this section is trying to capture the asymptotic behavior of the mountain pass solution.
\bl\label{le:230914-3}
Under the assumptions of Theorem \ref{th:230821-1}, $\{(\bar u_{\nu}, \bar v_{\nu})\}$ is a family of mountain pass solutions given by Theorem \ref{th:230821-1}. The following hold:
\begin{itemize}
\item [(i)] If $\mu_1<\mu_2$, then there exists $\varepsilon_{\nu}>0$ such that $(\bar u_{\nu},  \varepsilon_{\nu}^{\frac{N-2}{2}}\bar v_{\nu}(\varepsilon_{\nu}x))\to(0,\mu_{2}^{\frac{N-2}{4}}U)$ in $D_{0}^{1,2}(\R^N)\times D_{0}^{1,2}(\R^N)$ as $\nu\to0$.
\item [(ii)] If $\mu_1>\mu_2$, then there exists $\varepsilon_{\nu}>0$ such that $(\varepsilon_{\nu}^{\frac{N-2}{2}}\bar u_{\nu}(\varepsilon_{\nu}x), \bar v_{\nu})\to(\mu_{1}^{\frac{N-2}{4}}U,0)$ in $D_{0}^{1,2}(\R^N)\times D_{0}^{1,2}(\R^N)$ as $\nu\to0$.
\item [(iii)] If $\mu_1=\mu_2$, then one of the statements in (i) or (ii) holds.
\end{itemize}
Here $U$ is the unique positive radial solution of \eqref{eq:230921-1}.
\el
\bp
Since $\nu\to0$ and $0<M_{\nu}(a,b)<\frac{1}{N}\min\{{\mu_1}^{\frac{2-N}{2}},{\mu_2}^{\frac{2-N}{2}}\}S^{\frac{N}{2}}$, we get that $\{(\bar u_{\nu},\bar v_{\nu})\}$ is bounded in $E$. Assume that $\|\nabla \bar u_{\nu}\|_2\to l_1,\|\nabla \bar v_{\nu}\|_2\to l_2$, $l_1\geq0,l_2\geq0$. If $l_1=l_2=0$, then we have $M_{\nu}(a,b)=I_{\nu}(\bar u_{\nu},\bar v_{\nu})\to0$, which contradicts with $M_{\nu}(a,b)\geq k_0>0$. By $P_{\nu}(\bar u_{\nu},\bar v_{\nu})=0$ and $\int_{\R^N}\bar u_{\nu}^{\alpha}\bar v_{\nu}^{\beta}\ud x<+\infty$, we get that
$$\|\nabla \bar u_{\nu}\|_2^2+\|\nabla \bar v_{\nu}\|_2^2-\mu_1\|\bar u_{\nu}\|_{2^*}^{2^*}-\mu_2\|\bar v_{\nu}\|_{2^*}^{2^*}=o_{\nu}(1).$$
Combining with critical Sobolev inequality, we conclude that
$$l_1^2+l_2^2\leq (\mu_1 l_{1}^{2^*}+\mu_2 l_2^{2^*})S^{-\frac{2^*}{2}}.$$
Hence
$$l_1^2+l_2^2\geq \min\{l_1^2+l_2^2:0<l_1^2+l_2^2\leq S^{-\frac{2^*}{2}}(\mu_1 l_{1}^{2^*}+\mu_2 l_{2}^{2^*})\}\geq\min\{{\mu_1}^{\frac{2-N}{2}},{\mu_2}^{\frac{2-N}{2}}\}S^{\frac{N}{2}}.$$
Moreover, if $l_1^2+l_2^2=\min\{{\mu_1}^{\frac{2-N}{2}},{\mu_2}^{\frac{2-N}{2}}\}S^{\frac{N}{2}}$ then
\beq\label{230921-2}
\begin{cases}
(l_1,l_2)=(0,\mu_{2}^{\frac{2-N}{4}}S^{\frac{N}{4}}),&\quad\quad~\hbox{if}~ \mu_1<\mu_2,\\
(l_1,l_2)=(\mu_{1}^{\frac{2-N}{4}}S^{\frac{N}{4}},0),&\quad\quad~\hbox{if}~ \mu_1>\mu_2,\\
(l_1,l_2)=(0,\mu^{\frac{2-N}{4}}S^{\frac{N}{4}}) ~\hbox{or}~ (\mu^{\frac{2-N}{4}}S^{\frac{N}{4}},0), &\quad\quad~\hbox{if}~ \mu_1=\mu_2=\mu.
\end{cases}
\eeq
On the other hand, by $P_{\nu}(\bar u_{\nu},\bar v_{\nu})=0$ and $m_{\nu}(a,b)\to0$, we have
\begin{align*}
&\frac{1}{N}\min\{{\mu_1}^{\frac{2-N}{2}},{\mu_2}^{\frac{2-N}{2}}\}S^{\frac{N}{2}}\geq\lim\limits_{\nu\to0}M_{\nu}(a,b)=\lim\limits_{\nu\to0}I_{\nu}(\bar u_{\nu},\bar v_{\nu})\\
=&\frac{1}{N}(l_1^2+l_2^2)\geq\frac{1}{N}\min\{{\mu_1}^{\frac{2-N}{2}},{\mu_2}^{\frac{2-N}{2}}\}S^{\frac{N}{2}},
\end{align*}
which implies that $$\lim\limits_{\nu\to0}M_{\nu}(a,b)=\frac{1}{N}\min\{{\mu_1}^{\frac{2-N}{2}},{\mu_2}^{\frac{2-N}{2}}\}S^{\frac{N}{2}}$$
and  $l_1^2+l_2^2=\min\{{\mu_1}^{\frac{2-N}{2}},{\mu_2}^{\frac{2-N}{2}}\}S^{\frac{N}{2}}$.

Assume that $\mu_1<\mu_2$, then by \eqref{230921-2} we have $\|\nabla \bar u_{\nu}\|_2\to0$, $\|\nabla \bar v_{\nu}\|_2^2\to\mu_2^{\frac{2-N}{2}}S^{\frac{N}{2}}$ and $\|\bar v_{\nu}\|_{2^*}^{2^*}\to\mu_2^{-\frac{N}{2}}S^{\frac{N}{2}}$. Hence $\tilde v_{\nu}:=\mu_{2}^{\frac{2-N}{4}}\bar v_{\nu}$ is a minimizing sequence of \eqref{230913-3}. By \cite[Theorem 1.41]{Willem1996}, we know that there exists $\varepsilon_{\nu}>0$ such that $\varepsilon_{\nu}^{\frac{N-2}{2}}\tilde v_{\nu}(\varepsilon_{\nu}x)\to U$ in $D_{0}^{1,2}(\R^N)$. Hence $(\bar u_{\nu},\varepsilon_{\nu}^{\frac{N-2}{2}} \bar v_{\nu}(\varepsilon_{\nu}x))\to(0,\mu_2^{\frac{N-2}{4}}U)$ in $D_{0}^{1,2}(\R^N)\times D_{0}^{1,2}(\R^N)$.

For the cases of $\mu_1>\mu_2$ and $\mu_1=\mu_2$, we can prove them similar to the case of $\mu_1<\mu_2$ and we omit the details.
\ep

\vskip 0.2in
\noindent
{\bf Proof of Theorem \ref{th:230913-1}:}
It follows from Lemma \ref{le:230914-1}, Lemma \ref{le:230914-2} and Lemma \ref{le:230914-3}.
\hfill$\Box$


\end{document}